\documentclass[12pt]{article}
\newcommand{\qed}{\nolinebreak \hfill{$\Box $} \par\vspace{0.5\parskip}}

\newcommand{\rea}{\hspace{2pt}\hbox{R\hspace{-13pt} I\hspace{5pt}}}

\newcommand{\calb}{{\cal{B}}}

\newcommand{\cali}{{\cal{I}}}

\usepackage{txfonts}
\usepackage{graphicx}
\usepackage{color}

\newtheorem{theorem}{Theorem}[section]

\newtheorem{corollary}[theorem]{Corollary}

\newtheorem{definition}[theorem]{Definition}

\begin{document}

\begin{center}
\noindent
{\bf\Large A novel characterization and new simple tests of multivariate independence using copulas}

\noindent
{Jos\'e M. Gonz\'alez-Barrios\footnote{Corresponding author. 
E-mail address: gonzaba@sigma.iimas.unam.mx.}, 
Eduardo Guti\'errez-Pe\~na, Juan D. Nieves and Ra\'ul Rueda}

\noindent
{\it Department of Probability and Statistics, IIMAS, 
Universidad Nacional Aut\'onoma de M\'exico, 
Circuito Escolar s/n, Ciudad Universitaria, 
04510 Ciudad de M\'exico, Mexico}
\end{center}

\begin{abstract}
The purpose of this paper is twofold. First, we provide a novel characterization of independence of random vectors based on the checkerboard approximation to a multivariate copula. Using this result, we then propose a new family of tests of multivariate independence for continuous random vectors.
The tests rely on estimating the checkerboard approximation by means of the sample copula recently introduced in \cite{GBHC} and improved in \cite{GBHA2}. Such estimators have nice properties, including a Glivenko-Cantelli-type theorem that guarantees almost-sure uniform convergence to the checkerboard approximation.
Each of our test statistics is defined in terms of one of a number of different metrics, including the supremum, total variation and Hellinger distances, as well as the Kullback-Leibler divergence.
All of these tests can be easily implemented since the corresponding test statistics can be efficiently simulated under any alternative hypothesis, even for moderate and large sample sizes in relatively large dimensions. Finally, we assess the performance of our tests by means of a simulation study and provide one real data example.

\vspace{2ex}


\vspace{2ex}


\end{abstract}

\vskip .5cm

\begin{section}{\bf Introduction}

Mutual independence of a collection of random variables is a common assumption in many statistical procedures. Consequently, test of multivariate independence are required in a wide range of applications. Several test have been proposed, mostly for bivariate observations, and in the last few years there has been renewed interest in this topic, with some recent contributions focusing on the multivariate case. See, for example, \cite{BRS}, \cite{FLPS}, and the references therein.

Consider a $d$-dimensional random vector, where $d\geq 2$.
For the case $d=2$, popular tests for independence include those discussed by \cite{H}, \cite{BKR}, \cite{GR} and \cite{BP}. In some of these cases, the null distributions of the test statistics have been improved using suitable approximations; see \cite{MW}.

For $d\geq 3$, the number of proposals is significantly smaller. One reason for this is that, even though some of the above statistics can be extended to the multivariate case, the corresponding null distributions may become difficult to evaluate. On the other hand, some of the tests for $d=2$ are based on statistics such as the Kendall's $\tau$, the Spearman's $\rho$ or the Pearson's correlation coefficient, which do not have natural extensions to the case $d\geq 3$; see \cite{BRS} and \cite{TRO}.

In this paper, we propose a new class of test statistics for the general case $d\geq 2$ which  are competitive in the case $d=2$, while for the case $d\geq 3$ are feasible and easily implemented even for relatively large sample sizes. To achieve this, we first provide a simple novel characterization of $\Pi_d$ (the independence copula in dimension $d\geq 2$) in terms of the checkerboard approximations of orders $2$ and $3$. The checkerboard of order $m$, here denoted by $C^{(m)}$, is a  multilinear approximation  of a true $d$-copula $C$, based on the uniform partition of ${\bf I}=[0,1]$ given by $\{0,1/m,2/m,\ldots ,(m-1)/m,1\}$; see  \cite{CMM}, \cite{DFS} and \cite{Mik}. As we shall see, the independence copula $\Pi_d$, satisfies 
$\Pi_d(u_1,\ldots ,u_d)=C^{(m)}(u_1,\ldots ,u_d)$ for every $(u_1,\ldots ,u_d)\in {\bf I}^d$ and $m\geq 2$. However, for the converse we only need the equality to hold for $m=2$ and $m=3$.

The use of a sample $d$-copula of order $m$ (denoted $C_{n}^{(m)}$) as an estimator of $C^{(m)}$ was studied in \cite{GBHA2}.
Once the order $m$ has been chosen, the sample $d$-copula of order $m$ can be thought of  as a kind of kernel-based estimator of $C$. It turns out to be a good estimator of the copula $C$, even for moderate values of $m$, and improves on the empirical copula, the Bernstein copula and the beta empirical copula \cite{GBHA2}. On the other hand, the sample copula can be easily computed even for moderate sample sizes.

Let $H_0: C=\Pi_d$ denote the null hypothesis of independence;
that is, the true copula is the product copula. Since
$C_{n}^{(2)}$ and $C_{n}^{(3)}$ are unbiased estimators of $C^{(2)}$ and $C^{(3)}$, respectively, and
using the fact that $C_{n}^{(m)}$ and $C^{(m)}$ have constant densities on
the boxes of the uniform partition, we can propose
a test based on the distances between $C^{(2)}$ and $C_{n}^{(2)}$, and $C^{(3)}$ and $C_{n}^{(3)}$.
We consider several different distances including the supremum distance, the total variation distance and the Hellinger distance, as well as the Kullback-Leibler divergence. Our proposal works well for dimensions beyond $d = 2$. Moreover, the exact distributions of the test statistics can be  easily approximated using Monte Carlos methods since the sample copula is easy to simulate from. We shall show through simulations that our proposals have competitive powers in dimension $d=2$ and good powers when $d=3$ and $d=4$. We note, however, that the tests can be easily extended to higher dimensions.

The outline of the paper is as follows. In the next section, we review some necessary concepts relating to copulas and their estimation, including empirical copulas, checkerboard approximations, and sample copulas. This section also contains the main result of the paper, namely a new characterization of the independence copula based only on the checkerboards of orders $m=2$ and $m=3$. In Section~3, we describe in some detail the family of independence tests for copulas based on this characterization. Sections~4 and~5 present a simulation study and an application to a real dataset, respectively. Finally, in Section~6 we briefly discuss how our tests can be used to explore the dependence structure of random vectors.
\end{section}

\begin{section}{\bf A novel characterization of independence}
	
\begin{subsection}{\bf Preliminary results}

We start this section by reviewing some basic notions.

\begin{definition} \label{deftrivial}
Let ${\bf I}=[0,1]$ be the closed unit interval  and let $d\geq 2$ be an integer denoting the dimension.
Let $S_1,S_2,\ldots ,S_d$ be subsets of ${\bf I}$ such that $0,1\in S_i$ for every $i\in I_d$, where
$I_d=\{1,2,\ldots ,d\}$. The function $C:S_1\times S_2\times \cdots\times S_d\rightarrow \rea$ 
is a $d$-{\bf subcopula} if and only if $C$ satisfies:

i) $C(u_1,\ldots,u_d)=0$  if at least one $u_i=0$ for some $i\in I_d$.

ii) $C(1,\ldots,1,u_i,1,\ldots ,1)=u_i$ for every $i\in I_d$ and for every $u_i\in S_i$.

iii) $C$ is $d$-increasing, that is, for every $0\leq u_i\leq v_i\leq 1$ such that $u_i,v_i\in S_i$ for every $i\in I_d$, we have that if $B=[u_1,v_1]\times\cdots\times[u_d,v_d]$ then
\begin{equation}\label{volumen} 
V_C(B):= \sum_{\underline{b}} {\rm sgn}(\underline{b}) \, C(\underline{b})\geq 0,
\end{equation}
where the sum runs over all $\underline{b}=(b_1,\ldots,b_d)$ which are vertices of $B$, and the sign function is defined by
$$
{\rm sgn}(\underline{b})=\left\{\begin{array}{lcl}
1  & \mbox{if} & b_k=u_k \,\,\mbox{for an even number of} \,\,k's \\ \nonumber
-1 & \mbox{if} & b_k=u_k \,\,\mbox{for an odd number of} \,\, k's.
\end{array}\right.
$$
We say that a $d$-subcopula $C$ is a $d$-{\bf copula} if and only if $S_1=S_2=\cdots =S_d={\bf I}$.
\end{definition}

Recall that if $C$ is a $d$-copula then it is bounded below and above
by the  Fr\'echet-Hoeffding  bounds; that is, for every $(u_1,\ldots, u_d)\in{\bf I}^d$, we have
\begin{equation}\label{FreHoeff}
W_d(u_1,\ldots,u_d)\leq C(u_1,\ldots ,u_d)\leq M_d(u_1,\ldots ,u_d),
\end{equation}
where $W_d(u_1,\ldots,u_d)=\max\{ 0, u_1+\cdots+u_d-(d-1)\}$ and 
$ M_d(u_1,\ldots ,u_d)=\min\{u_1,\ldots ,u_d\}$. $M_d$ is always a $d$-copula, but
$W_d$ is only a copula for $d=2$. However, the left-hand inequality in (\ref{FreHoeff}) is always sharp.
Recall also that the {\bf independence copula} is defined by
$\Pi_d(u_1,\ldots,u_d)=u_1\cdot u_2\cdots u_d$ for every $(u_1,\ldots, u_d)\in{\bf I}^d$; see for example \cite{N06}.

Let $m$ be a positive integer, $I_m=\{1,2,\ldots ,m\}$, and define the {\bf uniform partition of size} $m$ {\bf of} ${\bf I}^d=[0,1]^d$ by
\begin{equation}\label{partunif}
R_{i_1,i_2,\ldots,i_d}^{m} = 
\left\langle \frac{i_1-1}{m},\frac{i_1}{m}\right] \times \left\langle \frac{i_2-1}{m},\frac{i_2}{m}\right]
\times\cdots\times \left\langle \frac{i_d-1}{m},\frac{i_d}{m}\right]
\,\,\mbox{for every}\,\, i_1,\ldots, i_d\in I_m,
\end{equation}
where the notation ``$\langle$'' means that
$\langle(i-1)/m,i/m]$ is a left open interval if $i>1$ or a closed interval if $i=1$.
Therefore, $\{ R_{i_1,\ldots ,i_d}^m\}_{(i_1,\ldots,i_d)\in I_m^d}$ is the uniform partition of ${\bf I}^d$ of size $m$.
Now let $C$ be a $d$-copula,  $m$ be a positive integer,
 $S_m=\{0,1/m,\ldots ,(m-1)/m, 1\}$ and define $C_m^{'}: S_m\times S_m \times\cdots \times S_m\rightarrow {\bf I}$ by
\begin{equation}\label{cm}
C_m^{'}=\{ C(u_1,\ldots ,u_d)\,\,|\,\, u_1,\ldots, u_d \in S_m\}.
\end{equation}
Clearly $C_m^{'}$ is $d$-subcopula.  So, if we define
$$ V_{C_m^{'}}(\overline{R_{i_1,\ldots ,i_d}^m}) =
\sum_{\underline{b}} \mbox{sgn}(\underline{b}) C(\underline{b}),
$$
where $\underline{b}$ denotes any vertex of $\overline{R_{i_1,\ldots ,i_d}^m}$, 
the closure of $R_{i_1,\ldots ,i_d}^m$,
then $V_{C_m^{'}}(\overline{R_{i_1,\ldots,i_d}^m})\geq 0$ for every $(i_1,\ldots, i_d)\in I_m^d$.

Using the multivariate version of Lemma 2.3.5 of Nelsen's book in the proof of Sklar's Theorem, see \cite{N06},  we know that the $d$-subcopula $C_m^{'}$ in equation (\ref{cm}) can be extended to a $d$-copula using a $d$-multilinear interpolation.
We denote this extended copula by $C^{(m)}$ and call it the 
{\bf checkerboard approximation of order $m$ of the $d$-copula $C$}.
It is well known, see \cite{CMM}, \cite{DFS}, \cite{LMT} and \cite{Mik}, that $C^{(m)}$ is a good approximation
of the true copula $C$ even for moderate  values of $m$. Moreover,
when $m\uparrow\infty$ we have a Glivenko-Cantelli-type theorem for the
supremum distance between $C^{(m)}$ and $C$.
In fact, $C^{(m)}$ is a good approximation to the true copula $C$. It is also easy to see that $C^{(m)}$ has a density function given by
\begin{equation}\label{denschecker}
c^{(m)}(u_1,\ldots ,u_d)=V_{C}(\overline{R_{i_1,\ldots ,i_d}^m})/ 
\lambda^d(\overline{R_{i_1,\ldots ,i_d}^m})=
m^d\cdot V_{C}(\overline{R_{i_1,\ldots ,i_d}^m})\,\,\mbox{for every}\,\,
(u_1,\ldots ,u_d)\in R_{i_1,\ldots ,i_d}^m,
\end{equation}
where $\lambda^d$ is the Lebesgue measure on the Borel space $({\bf I}^d,\calb({\bf I}^d))$.
Hence, the density is constant on each of the $d$-boxes of the uniform partition
for every $(i_1,\ldots i_d)\in I_{m}^{d}$.

We now provide an efficient method to estimate the checkerboard approximation of order $m$, that we call the {\bf sample $d$-copula of order $m$}. This methodology first appeared in \cite{GBHC} and was significantly improved in \cite{GBHA2}.
First, let $\underline{X}_{\,1}, \ldots ,\underline{X}_{\,n}$ be a random sample
from a continuous \mbox{$d$-dimensional} distribution function $H$ or a $d$-copula $C$; then the {\bf modified sample} or {\bf pseudo-sample} is defined by transforming the coordinates of the sample using ranks, and dividing
them by the sample size $n$. We denote the modified sample by
$U_n=\{ \underline{Y}_{\,1}, \ldots ,\underline{Y}_{\,n}\}$, where 
$\underline{Y}_{\,i}=(Y_{i,1},\ldots ,Y_{i,d})$
for every $i\in\{1,\ldots ,n\}$. Of course, $U_n$ is always a sample in ${\bf I}^d$.
We define the {\bf empirical copula (subcopula)} denoted by $C_n: {\bf I}^d\rightarrow {\bf I}$ by
\begin{equation}\label{empirica}
C_n(u_1,\ldots ,u_d)=\frac{1}{n}\sum_{i=1}^{n}{\bf 1}_{\{Y_{i,1}\leq u_1,\ldots ,Y_{i,d}\leq u_d\}}
(u_1,\ldots ,u_d)\,\,\mbox{for every}\,\,(u_1,\ldots ,u_d)\in {\bf I}^d,
\end{equation}
where ${\bf 1}_{A}$ denotes the indicator function of the set $A$.
This is actually a subcopula on the grid $\{0,1/n,\ldots (n-1)/n,1\}^d$, because it has jumps of magnitude $1/n$ at every modified observation $\overline{Y}_{\,i}$. It is well known that the empirical copula converges almost surely and uniformly  to the true copula $C$ as $n$  increases to infinity, see for example \cite{N06}. 

\begin{definition}\label{samplecopula}
Let $U_n=\{ \underline{Y}_{\,1}, \ldots ,\underline{Y}_{\,n}\}$ be a modified sample
of size $n$ from a joint $d$-dimensional distribution function $H$ or a $d$-copula $C$.
Let  $2\leq m\leq n$ be an integer, such that $m$ divides $n$, and let
$R_{i_1,i_2,\ldots,i_d}^{m} = 
\left\langle \frac{i_1-1}{m},\frac{i_1}{m}\right] \times \left\langle \frac{i_2-1}{m},\frac{i_2}{m}\right]
\times\cdots\times \left\langle \frac{i_d-1}{m},\frac{i_d}{m}\right]$ be a box of the uniform
partition of size $m$ from ${\bf I}^d$. Define, for every $\underline{i}=(i_1,\ldots, i_d)\in I_{m}^{d}$,
\begin{equation}\label{eses}
s_{i_1,\ldots ,i_d}^{n,(m)}=\frac{{\rm card}(R_{\underline{i}}^{m}\cap U_n)}{n},
\end{equation}
where ${\rm card}(\cdot)$ denotes cardinality of a set, that is, $s_{i_1,\ldots ,i_d}^{n,(m)}$
are the relative frequencies in each $d$-box. Finally, let
\begin{equation}\label{matrizS}
S_{m}^{n} =\left(s_{i_1,\ldots ,i_d}^{n,(m)}\right)_{(i_1,\ldots ,i_d)\in I_{m}^{d}}.
\end{equation}
Then $S_{m}^{n}$ is a $d$-dimensional matrix whose entries add up to $1$. If we
define, for every $j\in \{1,2,\ldots ,d\}$ and $k\in I_m$,
$$ q_{j,k}=\sum_{i_j=1}^{k}\sum_{i_1=1}^{m}\cdots \sum_{i_{j-1}=1}^{m}
\sum_{i_{j+1}=1}^{m} \cdots \sum_{i_d=1}^{m} s_{i_1,\ldots ,i_d}^{n,(m)},$$
then $0=q_{j,0}<\frac{1}{m}=q_{j,1}<\cdots <\frac{m-1}{m}=q_{j,m-1}<1=q_{j,m}$ 
for every $j\in \{1,2,\ldots ,d\}$; that is, we recover the
uniform partition of order $m$ of \ ${\bf I}^d$. Now let 
\begin{equation}\label{subcopulaS}
t_{i_1,\ldots ,i_d}=\sum_{k_1=1}^{i_1}\sum_{k_2=1}^{i_2}\cdots \sum_{k_{d-1}=1}^{i_{d-1}}
\sum_{k_d=1}^{i_d} s_{k_1,\ldots ,k_d}^{n,(m)},
\end{equation}
for every point $(i_1,\ldots ,i_d)$ on the grid $I_m^d$; 
then $T =\left(t_{i_1,\ldots ,i_d}\right)_{(i_1,\ldots ,i_d)\in I_{m}^{d}}$ is a $d$-subcopula on the grid $\{0,1/m,\ldots ,(m-1)/m,1\}^d$.
Here it is convenient to recall the $d$-multilinear interpolation in the proof of Sklar's Theorem for $d=2$ (see equation (2.3.2) of Lemma 2.3.5 in Nelsen's book \cite{N06}). If we apply this interpolation to the $d$-subcopula $T$, the result is what we call the {\bf sample} $d$-{\bf copula of order} $m$. It will be denoted by $C_{n}^{(m)}$.
\end{definition}

It can be shown, see \cite{GBHA2}, that $C_{n}^{(m)}$ is always a $d$-copula and that it estimates $C^{(m)}$. 
In fact, $C_{n}^{(m)}$ is related to a type of kernel estimator of the true copula $C$ when the density of $C$ exists.
Moreover, $C_{n}^{(m)}$  always has a density $c_{n}^{(m)}(u_1,\ldots ,u_d)$ 
which is constant for every $(i_1,\ldots, i_d)\in I_{m}^{d}$ and is given by

\begin{equation}\label{densitycmn}
c_{n}^{(m)}(u_1,\ldots ,u_d) = \frac{s_{i_1,\ldots ,i_d}^{n,(m)}}
{\lambda^d(\overline{R_{i_1,\ldots ,i_d}^m})}= m^d\cdot s_{i_1,\ldots ,i_d}^{n,(m)}
\,\mbox{for every}\,(u_1,\ldots ,u_d)\in R_{i_1,\ldots ,i_d}^m.
\end{equation}

It is worth pointing out that the sample $d$-copula of order $m$ is far easier to compute than the empirical copula, the Bernstein copulas and the beta empirical copulas, see \cite{JSV} and \cite{SST}, all of which  have been used to estimate the true copula $C$. However, as shown in \cite{GBHA2}, in all of these
cases we can obtain better approximations to the true copula $C$ using the sample copula.

We also have a Glivenko-Cantelli-type theorem which gives uniform almost sure convergence
of $C_{n}^{(m)}$ to $C^{(m)}$, for every $m\geq 2$, that is,
\begin{equation}\label{GCm}
\lim_{n\rightarrow\infty} \sup_{(u,v)\in{\bf I}^2} |C_{n}^{(m)}(u,v)-C^{(m)}(u,v)|=0\,\,\mbox{a.s.}
\end{equation}
On the other hand, from \cite{CMM}, we also have
\begin{equation}\label{ConvC^m}
\lim_{m\rightarrow\infty} \sup_{(u,v)\in{\bf I}^2} |C^{(m)}(u,v)-C(u,v)|\leq 
\lim_{m\rightarrow\infty} \frac{d}{2m} =0.
\end{equation}

Now let $P_{C_{n}^{(m)}}$ and $Q_{C^{(m)}}$ be the probability measures induced by the
sample $d$-copula $C_{n}^{(m)}$ and by the checkerboard copula $C^{(m)}$, respectively, associated
to a $d$-copula $C$. Recall that the {\bf total variation distance}, see for example
\cite{GS}, between two probability measures
$P$ and $Q$ on the Borel measurable space $({\bf I}^d,\calb({\bf I}^d))$
is defined by

\begin{equation}\label{tv}
d_{\mbox{\tiny{\it TV}}}(P,Q) =\sup_{A\in\calb({\bf I}^d)} |P(A)-Q(A)|\leq 1.
\end{equation}
Suppose that $P$ and $Q$ have densities $f_P$ and $f_Q$, with respect to the Lebesgue measure $\lambda^d$ on the measurable space $({\bf I}^d, \calb({\bf I}^d))$, which are constant on the uniform partition of order $m$ of ${\bf I}^d$. Then the total variation distance of $P$ and $Q$ can be written as
\begin{equation}\label{tvdiscrete}
d_{\mbox{\tiny{\it TV}}}(P,Q) =\sup_{A\in\calb({\bf I}^d)} |P(A)-Q(A)|=\frac{1}{2} \int_{{\bf I}^d}
|f_P-f_Q|d\lambda^d.
\end{equation}
Using equations (\ref{denschecker}), (\ref{densitycmn}), (\ref{GCm}) 
together with equation (\ref{tvdiscrete}), it is easy to prove the following.

\begin{theorem}\label{TGCdTV}
Let $C$ be a $d$-copula. Take $m\geq 2$ fixed and let $n$ be a multiple of $m$.
Denote by $C_{n}^{(m)}$ the sample copula of order $m$ built from a modified sample of size $n$ from $C$, and let $C^{(m)}$ be the corresponding checkerboard approximation of order $m$. If $P_{C_{n}^{(m)}}$ and $Q_{C^{(m)}}$ are the 
probability measures on ${\bf I}^d$ defined by $C_{n}^{(m)}$ and $C^{(m)}$, respectively, then
\begin{equation}\label{GCdTV}
\lim_{n\rightarrow \infty} d_{\mbox{\tiny{\it TV}}}\left(P_{C_{n}^{(m)}},Q_{C^{(m)}}\right)=0\quad \mbox{a.s.}
\end{equation}
\end{theorem}

Other important metrics are the {\bf Hellinger distance} and the 
{\bf supremum distance} or {\bf uniform distance}, see \cite{GS}. 
The first one is a $L^2$-type distance between $P$ and $Q$; it is defined in terms of the corresponding density functions $f_P$ and $f_Q$, and is given by 
\begin{equation}\label{Helli}
d_H(P,Q) = \frac{1}{\sqrt{2}}
\left[\int_{{\bf I}^d}\left(\sqrt{f_P}-\sqrt{f_Q}\right)^2 d\lambda^d\right]^{1/2}\leq 1.
\end{equation}
The second one is also called  the {\bf weak distance}, because it is related to
weak convergence. Let $F_P$ and $F_Q$ be the distribution functions associated
to the probability measures $P$ and $Q$, respectively Then
\begin{equation}\label{supdist}
d_{\mbox{sup}}(F_P,F_Q)=d_{\infty}(F_P,F_Q)=\sup_{\underline{x}\in {\bf I}^d}
\left|F_P(\underline{x})-F_Q(\underline{x})\right|\leq 1.
\end{equation}

Finally, we also consider one more functional that is not a metric, namely the {\bf relative entropy},
also known as the {\bf Kullback-Leibler divergence}. For two probability measures
$P$ and $Q$ with densities $f_P$ and $f_Q$, this is given by
\begin{equation}\label{KL}
d_{I}(P,Q)=\int_{S(P)} f_P\log\left(\frac{f_P}{f_Q}\right) d\lambda^d,
\end{equation}
where $S(P)$ is the support of $P$ on $\rea^d$, and we define $0\log(0/q)=0$ for every $q\in\rea$
and $p\log(p/0)=\infty$, see \cite{GS}.
This divergence satisfies $d_I(P,P)=0$ and $d_{I}(P,Q)\geq 0$, but 
it is not symmetric and does not satisfy the triangle inequality. However, it is an important quantity in Statistics, as it measures information loss.

In the next section, we shall use equations (\ref{tvdiscrete}), (\ref{Helli}), (\ref{supdist}) and (\ref{KL}) to define four statistics to test for multivariate independence.

\end{subsection}

\begin{subsection}{\bf Main result}

In this section we present a characterization of independence in terms of checkerboard approximations of a multivariate copula.

\begin{theorem} \label{MainThm}
Let $C$ be a d-copula. Then $C=\Pi_d$ if and only if
\begin{equation}\label{result}
C(u_1,\ldots,u_d)=
C^{(2)}(u_1,\ldots,u_d)=C^{(3)}(u_1,\ldots,u_d)\,\,\mbox{ for all }\,\,(u_1,\ldots,u_d)\in{\bf I}^d,
\end{equation}
where $C^{(2)}$ and $C^{(3)}$ are the checkerboard approximations of the $d$-copula $C$ of order 2 and 3, respectively.
\end{theorem}

\noindent The proof of this theorem is given in the Appendix.

\vspace{2ex}

From equations (\ref{gral}) and (\ref{gral3}) in the proof we have the following.

\begin{corollary} \label{coro} 
Let $C=\Pi_d$ be the product copula. Then, for $m\geq 2$, we have
$$ C^{(m)}(u_1,\ldots,u_d)=\Pi_d(u_1,\ldots ,u_d) \, \mbox{ for all } \, (u_1,\ldots ,u_d)\in{\bf I}^d.$$
\end{corollary}
\end{subsection}
\end{section}

\begin{section}{Independence tests}

The total variation  distance defined in equations (\ref{tv}) and (\ref{tvdiscrete})
provides the largest possible difference between two probability measures,
so it is considered a stronger distance than the ``sup'' distance.  The total variation distance is often regarded as ``too strong to be useful'', but this is not so in our case, as Theorem~\ref{TGCdTV} shows. 

Using the characterization of independence given in Theorem~\ref{MainThm}, we first propose
a  new independence test based on the total variation distance. We know by equation (\ref{result}) that,
for $d\geq 2$,
$$ C=\Pi_d \quad\mbox{if and only if}\quad C=C^{(2)}=C^{(3)}.$$
Let $Q_{C^{(2)}}$ and $Q_{C^{(3)}}$ be the probability measures
induced by the checkerboards of order $m=2$ and $m=3$, respectively. Assuming (\ref{result}) holds, if we observe that the probability measure associated with $\Pi_d$ is simply the Lebesgue product measure $\lambda^d$ then we have
\begin{equation}\label{dtvProd}
d_{\mbox{\tiny{\it TV}}}(Q_{C^{(2)}},\lambda^d)=0\quad \mbox{and}\quad
 d_{\mbox{\tiny{\it TV}}}(Q_{C^{(3)}},\lambda^d)=0.
\end{equation}

Since the total variation distance is quite strong, we may use it to see whether the
true copula $C$ equals the product copula $\Pi_d$ or not. Thus we shall use  
the fact that, under $H_0: C=\Pi_d$, the measures $Q_{C^{(2)}}$ and $Q_{C^{(3)}}$ 
are equal to $\lambda^d$ (by equation (\ref{dtvProd})). Besides, by
Theorem~\ref{TGCdTV}, $Q_{C^{(2)}}$ and $Q_{C^{(3)}}$ are the uniform limits
of $P_{C_{n}^{(2)}}$ and $P_{C_{n}^{(3)}}$ as $n$ increases. Hence, based on Corollary~\ref{coro}
we propose the statistic
\begin{equation}\label{stattvd}
\eta_{\mbox{\tiny{\it TV}}}(C;n) =
\frac{ d_{\mbox{\tiny{\it TV}}}(P_{C_{n}^{(2)}},\lambda^d)+
d_{\mbox{\tiny{\it TV}}}(P_{C_{n}^{(3)}},\lambda^d)}{2}.
\end{equation}
In this case, we take a sample of sample size $n$ from the true copula $C$. Now recall that $P_{C_{n}^{(m)}}$, for $m=2$ and $m=3$, are the probability measures induced by the sample copulas $C_{n}^{(m)}$ of orders $m=2$ and $m=3$, respectively. Note that, by equation (\ref{GCdTV}), $\lim_{n\rightarrow\infty} \eta_{\mbox{\tiny{\it TV}}}(C;n)=0$ almost surely under $H_0$.
On the other hand, the alternative hypothesis is $H_1: C\not = \Pi_d$ so for any copula $C\not = \Pi_d$ we have $\lim_{n\rightarrow\infty} \eta_{\mbox{\tiny{\it TV}}}(C;n)>0$.

Even though the null distribution of the test statistic $\eta_{\mbox{\tiny{\it TV}}}(\Pi_d;n)$ is generally not known for a fixed sample size $n$, it is straightforward to generate a large number of simulated samples (under $H_0$) even for moderate values of the dimension $d$ and  sample size $n$. We can then use these samples to approximate the quantiles of order $90\%, 95\%$ and $99\%$, say, needed to perform a standard test. We would then reject $H_0$ at level $\alpha$ if the
observed value of $\eta_{\mbox{\tiny{\it TV}}}(C;n)$ exceeds the corresponding $(1-\alpha)$ quantile.

Note that we can replace the total variation distance with any of the other measures discussed at the end of Section~2. For example, we can use the  Hellinger distance given in equation (\ref{Helli}), or the supremum distance given in equation (\ref{supdist}). Furthermore, we can even use the
Kullback-Leibler divergence of equation (\ref{KL}).  Since  the densities of the
product copulas $\Pi_d$, and the sample $d$-copulas of order $m$ are constants on the $d$-boxes of the uniform partition as in equation (\ref{partunif}),  $d_{I}$  satisfies 
$d_{I}(Q_{C^{(2)}},\lambda^d)=0=d_{I}(Q_{C^{(3)}},\lambda^d)$ if $C=\Pi_d$. On the other hand, both are strictly positive if $C\not= \Pi_d$. Therefore, we can also use any of the following  statistics to test
for multivariate independence.

\begin{equation}\label{statHel}
\eta_{\mbox{\tiny{\it H}}}(C;n) =
\frac{ d_{H}(P_{C_{n}^{(2)}},\lambda^d)+d_{H}(P_{C_{n}^{(3)}},\lambda^d)}{2},
\end{equation}

\begin{equation}\label{statsup}
\eta_{\mbox{\tiny{\it sup}}}(C;n) =
\frac{ d_{\mbox{sup}}(C_{n}^{(2)},\Pi_d)+d_{\mbox{sup}}(C_{n}^{(3)},\Pi_d)}{2},
\end{equation}

\begin{equation}\label{statKL}
\eta_{\mbox{\tiny{\it I}}}(C;n) =
\frac{ d_{I}(P_{C_{n}^{(2)}},\lambda^d)+d_{I}(P_{C_{n}^{(3)}},\lambda^d)}{2}.
\end{equation}

In \cite{GBHA1} both the exact and the asymptotic distributions of $C_{n}^{(m)}$ under independence are derived. However, finding the null distribution of the test statistics is not straightforward, 
which is why we resort to Monte Carlo tests.

\end{section}

\begin{figure}[t!]
\begin{center}
\includegraphics[width=5.1cm,height=5.1cm]{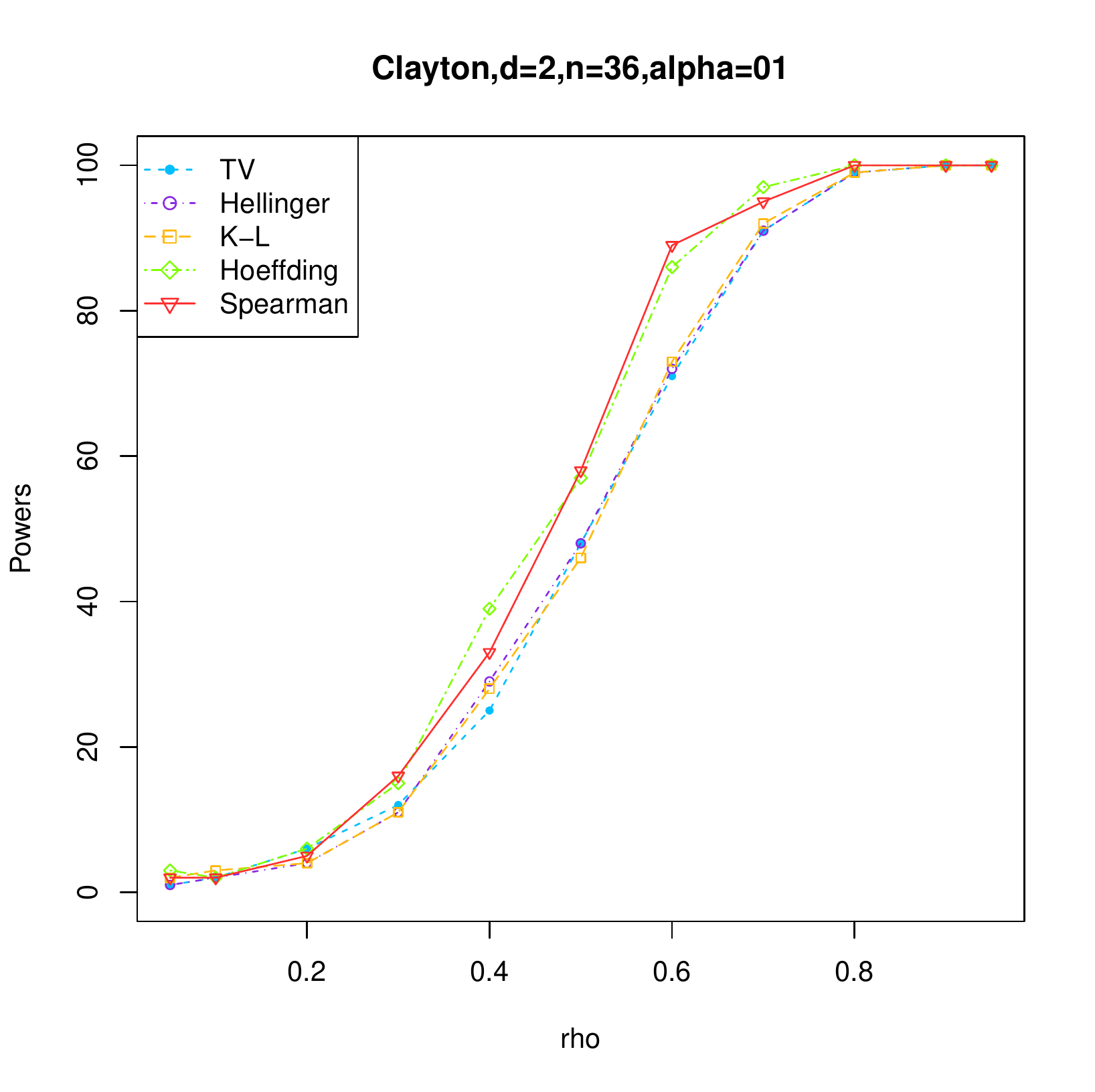}\hspace{.1cm}
\includegraphics[width=5.1cm,height=5.1cm]{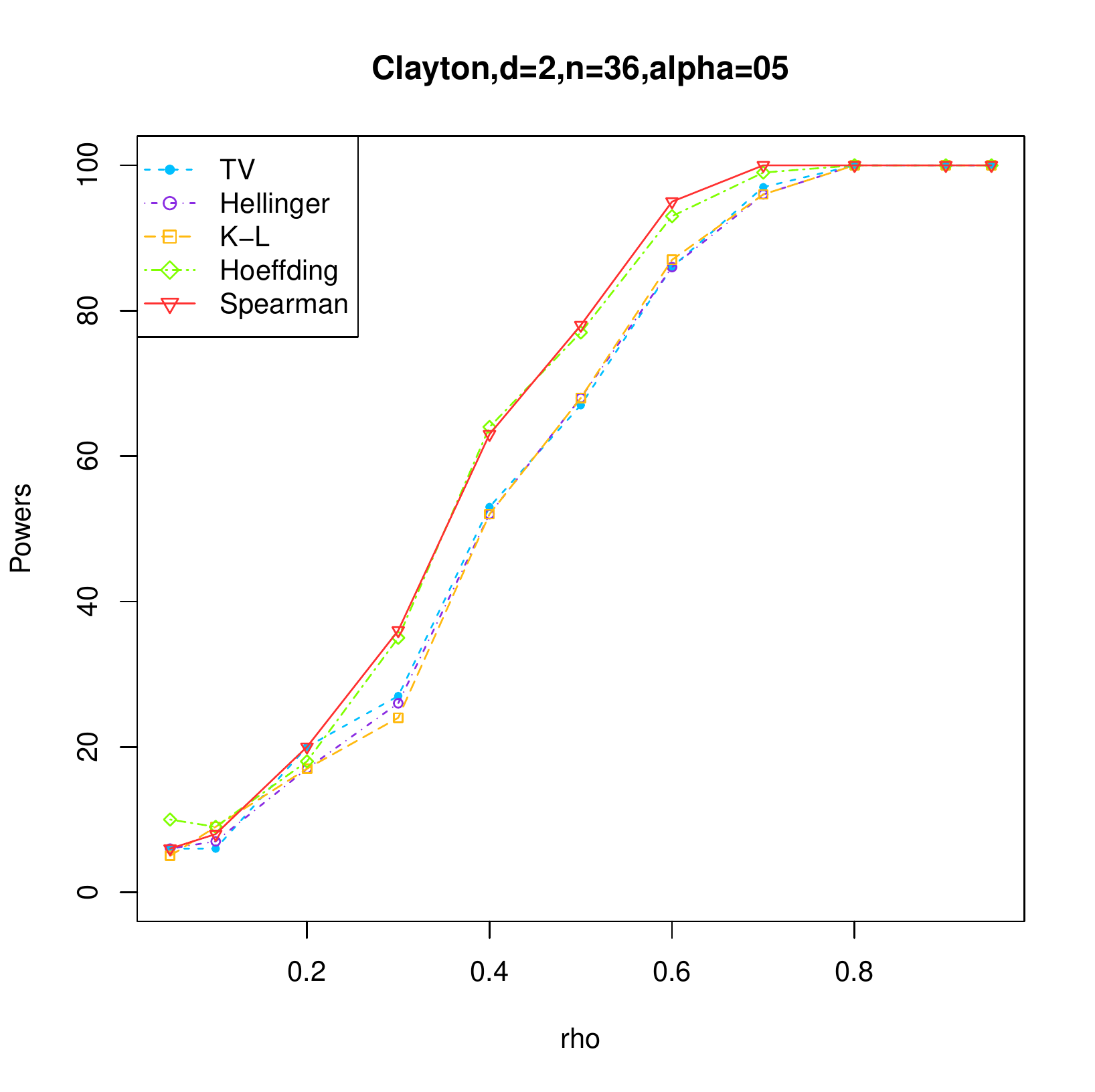}\hspace{.1cm}
\includegraphics[width=5.1cm,height=5.1cm]{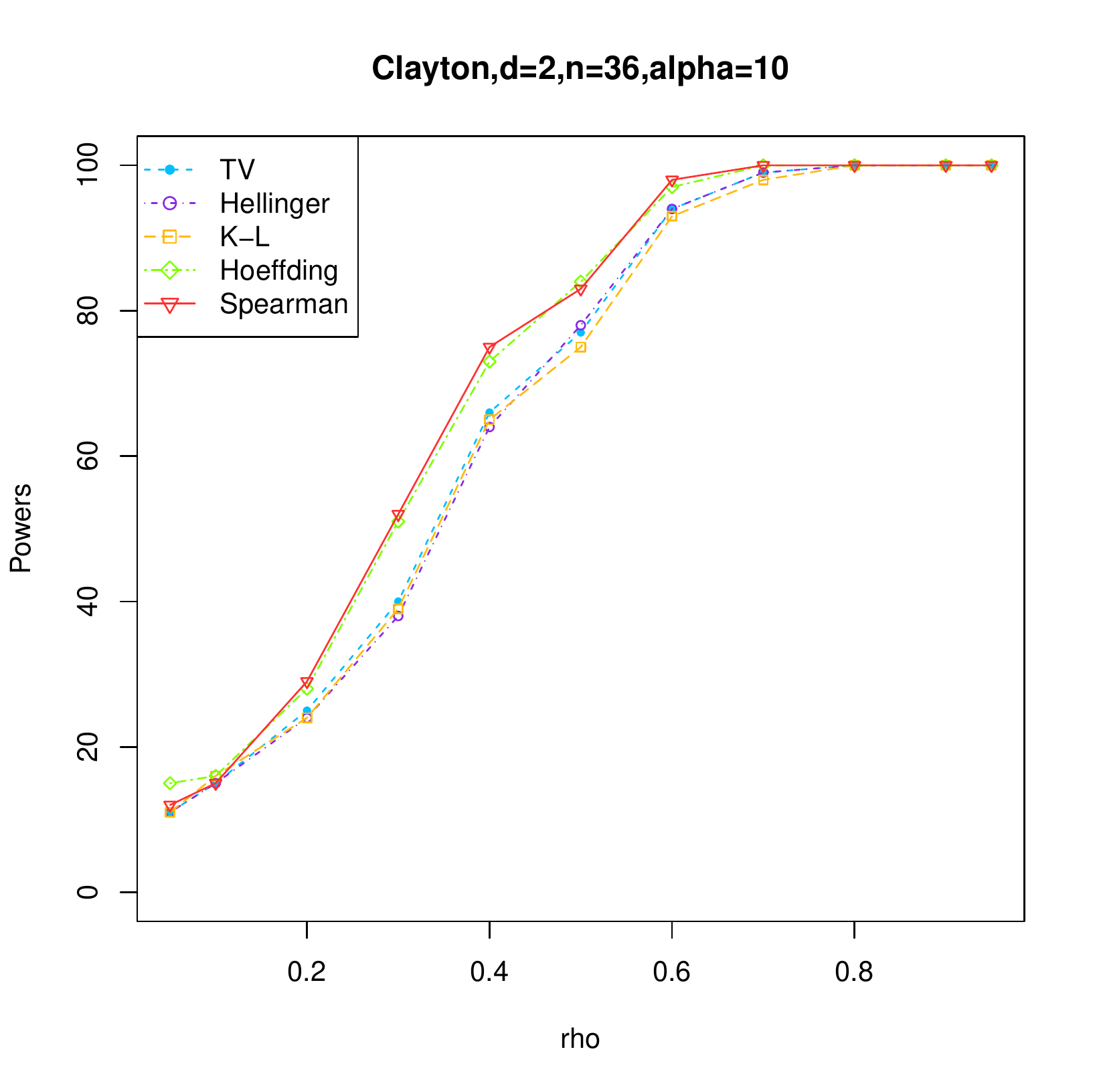}\hspace{.1cm}
\caption{Power comparisons for the Clayton family; $n=36$, $d=2$.}
\label{Clayton}
\end{center}
\end{figure}

\begin{figure}[h]
\begin{center}
\includegraphics[width=5.1cm,height=5.1cm]{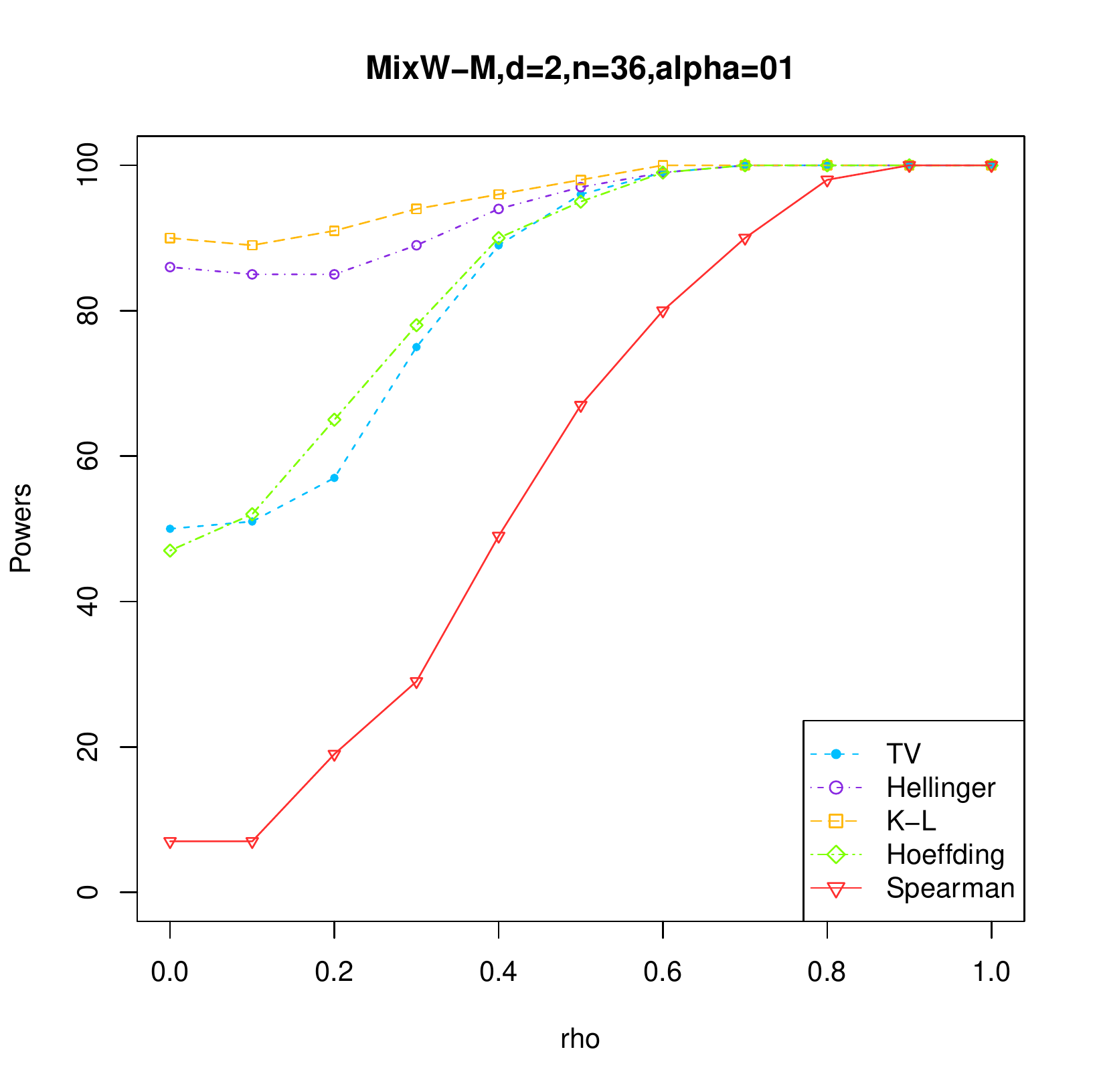}\hspace{.1cm}
\includegraphics[width=5.1cm,height=5.1cm]{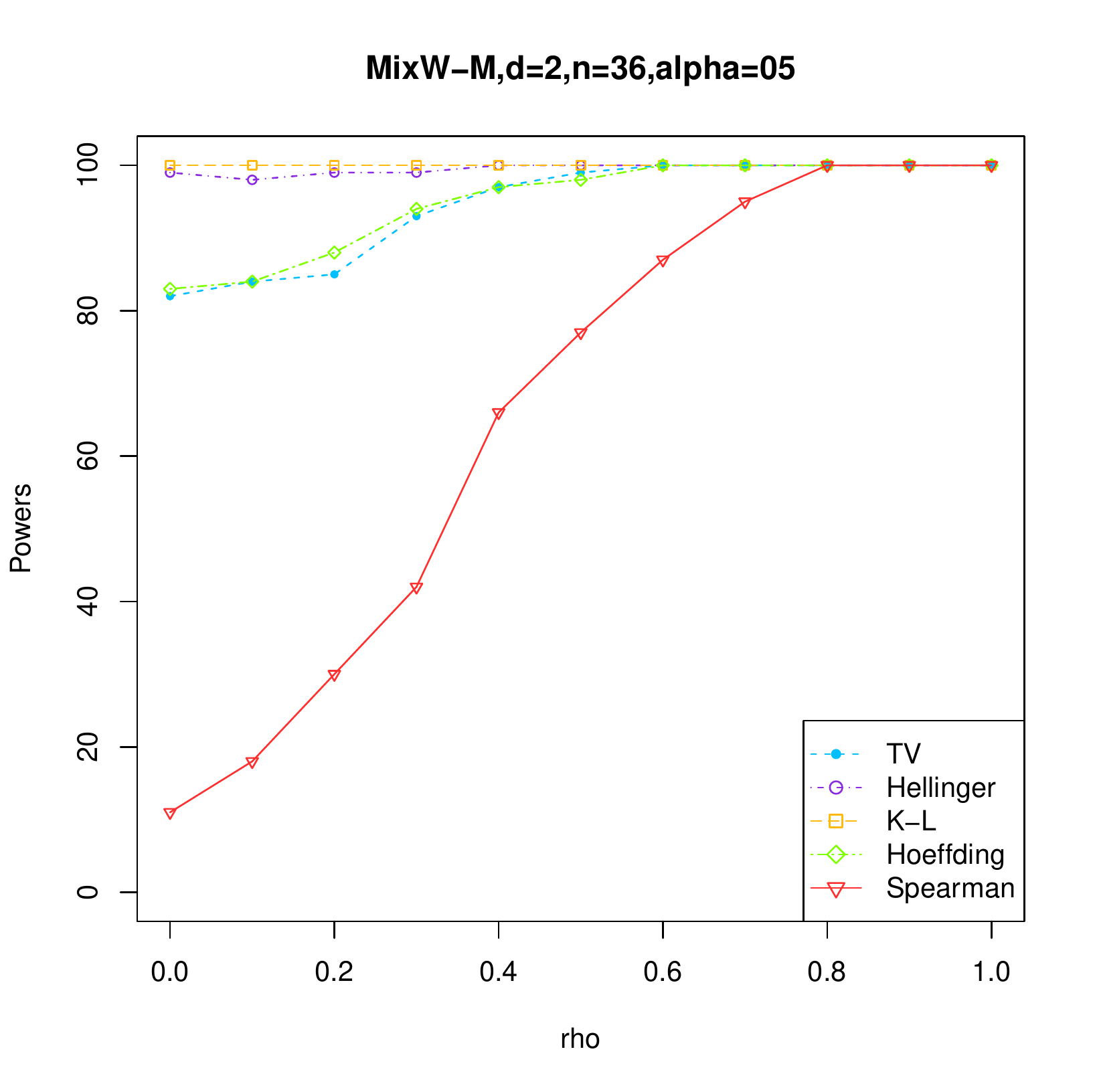}\hspace{.1cm}
\includegraphics[width=5.1cm,height=5.1cm]{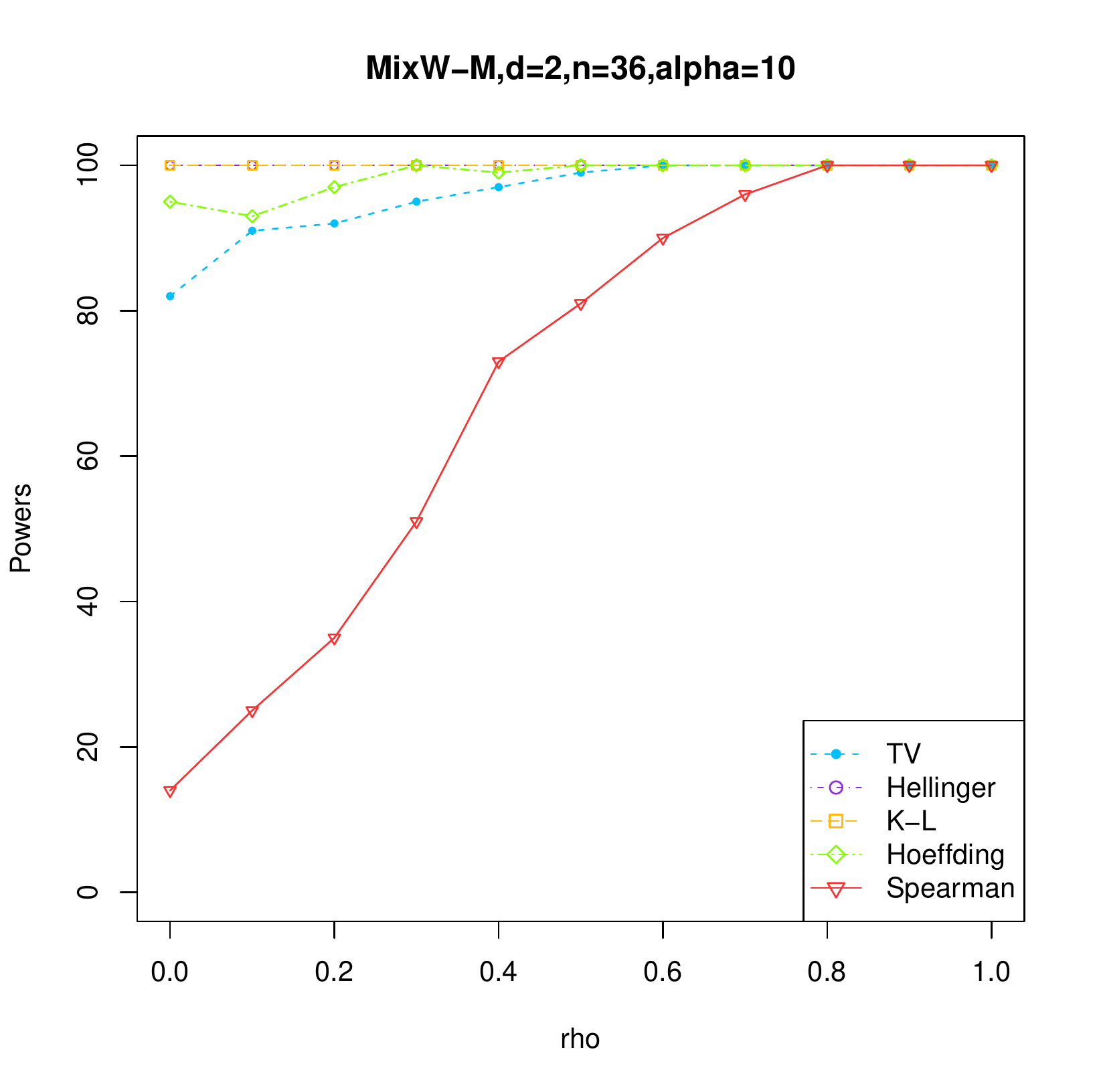}\hspace{.1cm}
\caption{Power comparisons for the Frechet-Mardia family; $n=36$, $d=2$.}
\label{Fre-Mar}
\end{center}
\end{figure}

\begin{figure}[h!]
\begin{center}
\includegraphics[width=5.1cm,height=5.1cm]{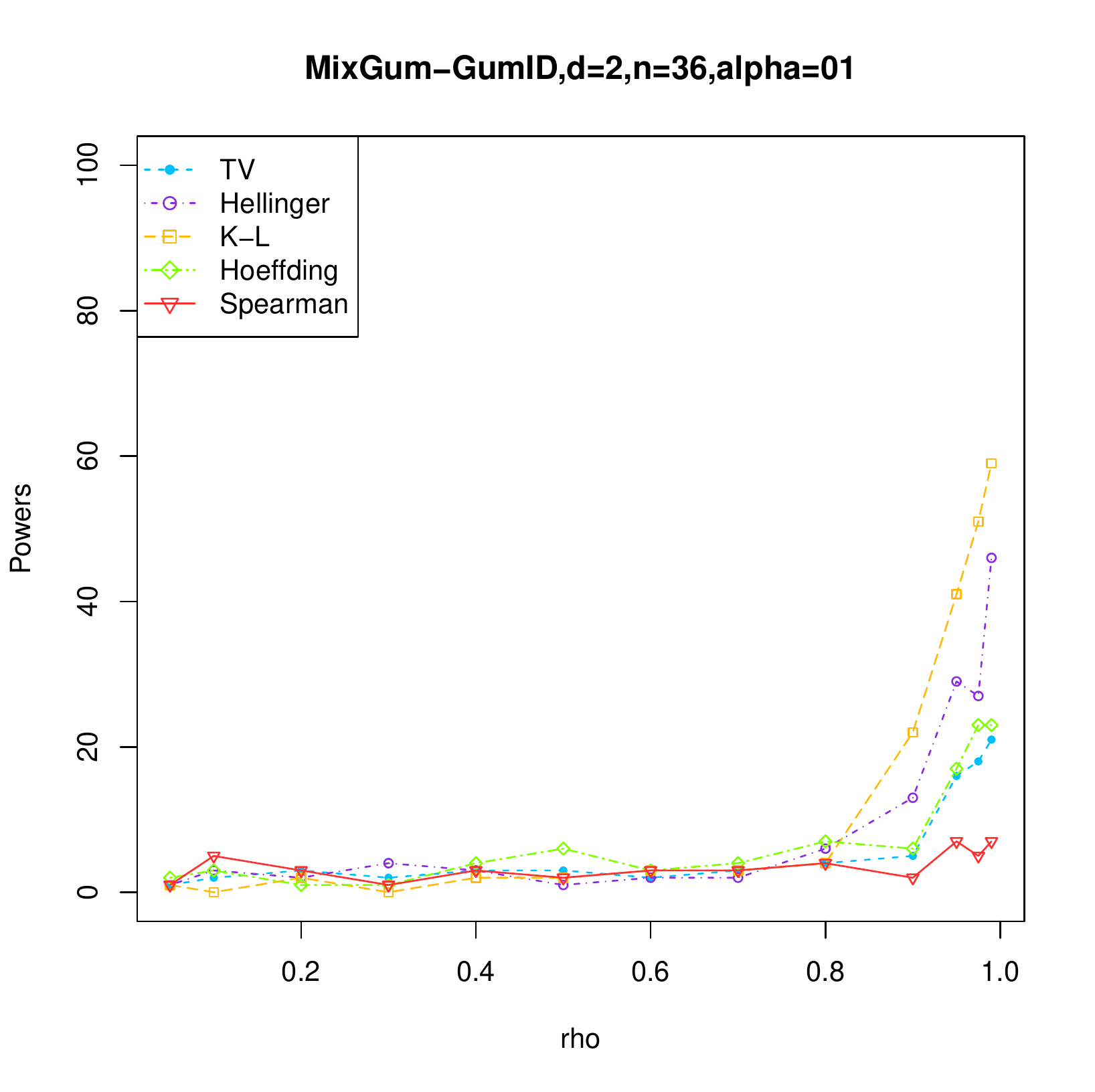}\hspace{.1cm}
\includegraphics[width=5.1cm,height=5.1cm]{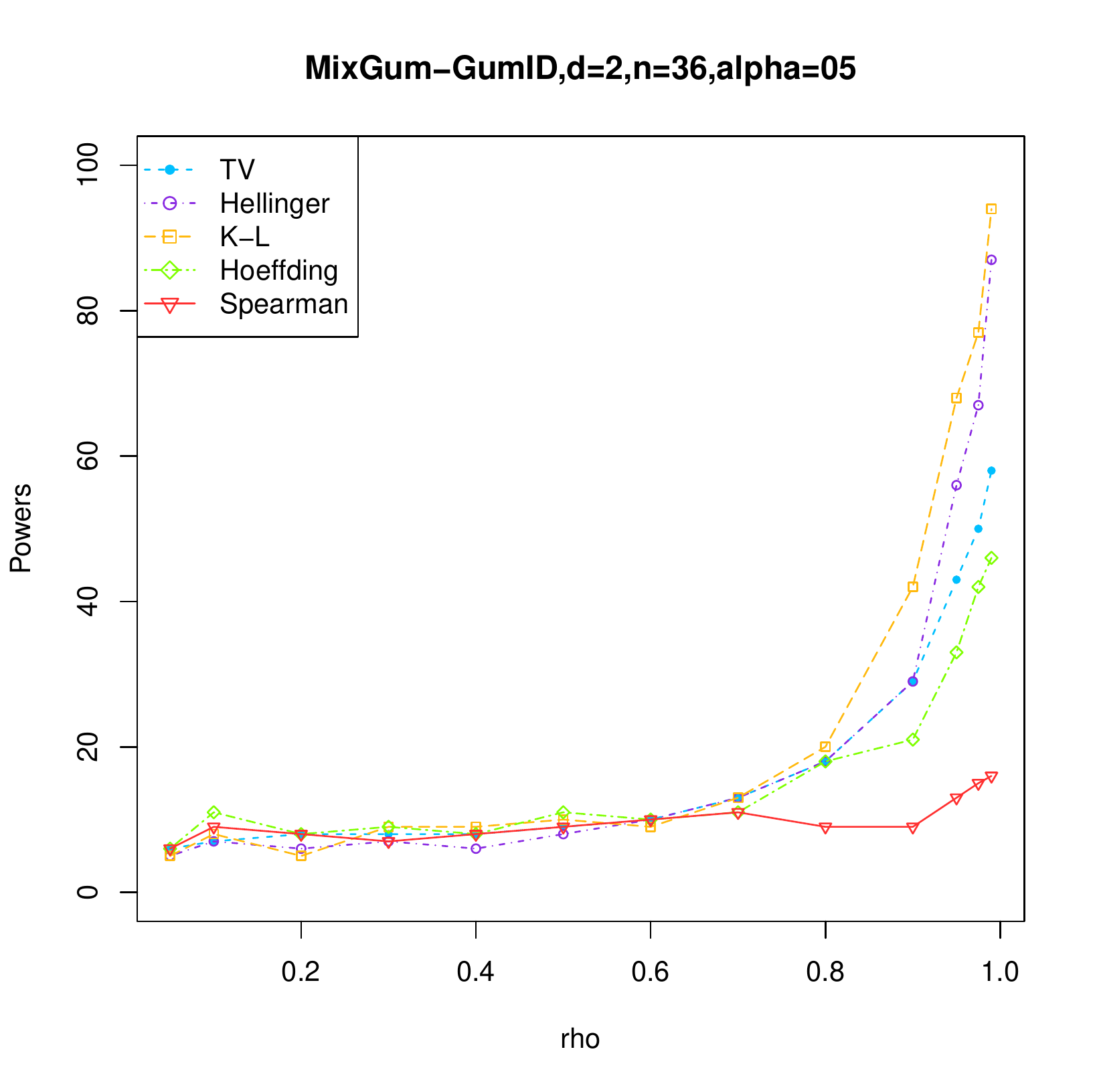}\hspace{.1cm}
\includegraphics[width=5.1cm,height=5.1cm]{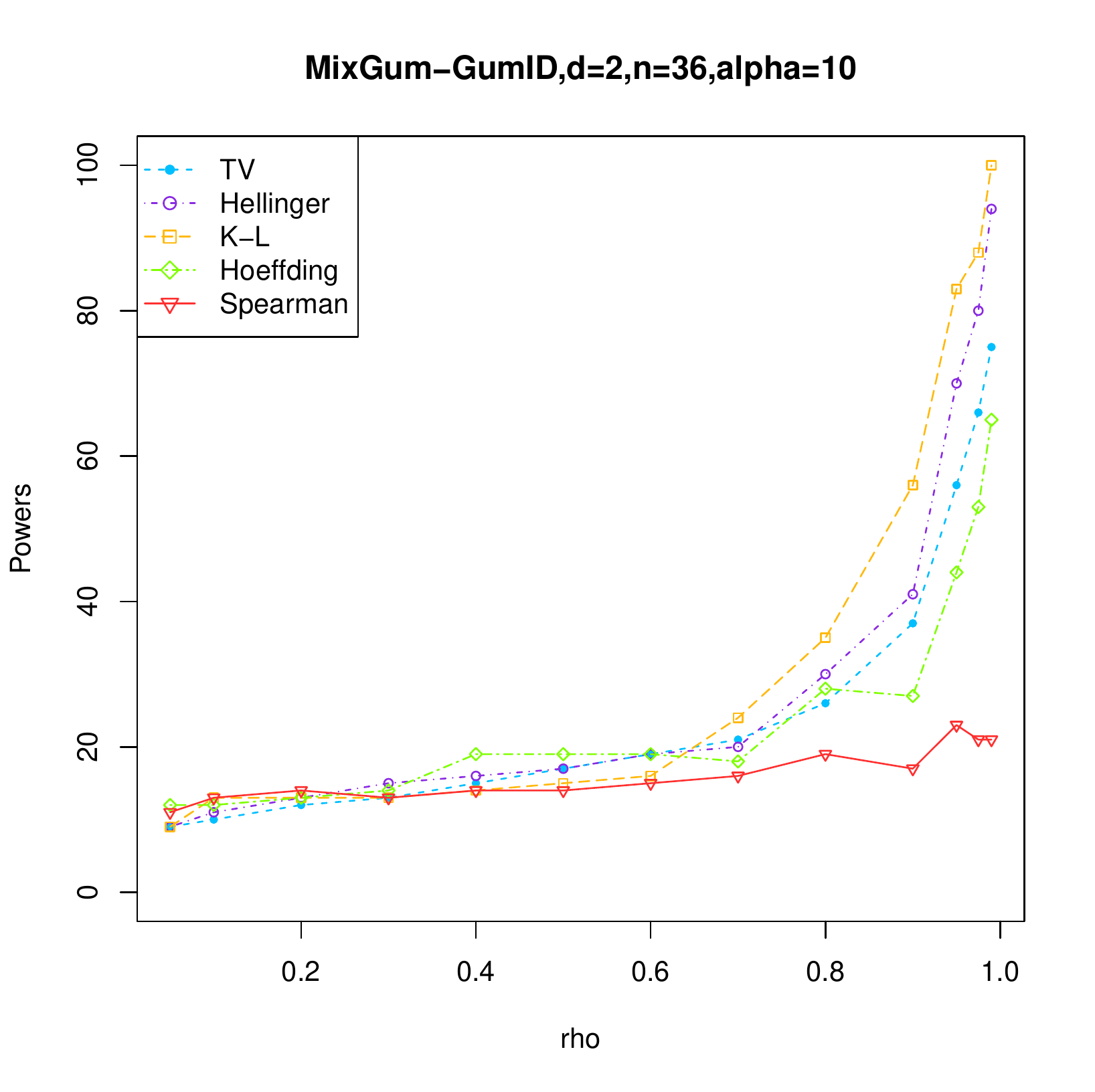}\hspace{.1cm}
\caption{Power comparisons for the Mixture Gumbel-Gumbel ID family; $n=36$, $d=2$.}
\label{Gum-GumID}
\end{center}
\end{figure}

\begin{section}{Simulation study}

In this section we will carry out a simulation study in dimensions $d=2$, $d=3$ and $d=4$. We start by comparing our tests with several well known tests of independence in the case $d=2$. Then, 
for dimensions $d=3$ and $d=4$ we present the results of a comparison among our test statistics (\ref{stattvd})--(\ref{statKL}) for several different families of copulas.

\begin{subsection}{Dimension $d=2$}

Several statistics have been proposed to test for independence between two random variables. Here, we consider two classical tests. The first one was proposed by Hoeffding in 1948 to test the independence of two continuous random variables with continuous joint and marginal densities, 
see \cite{H}. It is based on the function $D(x,y)=F(x,y)-F(x,\infty)\cdot F(\infty,y)= F(x,y)-F_1(x)\cdot F_2(y)$, where
$F$ denotes the joint distribution function, and $F_1$ and $F_2$
are the marginal distributions of $X$ and $Y$. 
The test statistic he proposed is based on the empirical version of 
$\Delta(F)=\int D^2(x,y) dF(x,y)$.
Here we used the {\bf hoeffd} function of the {\sc R package Hmisc} (R Core Team \cite{RCT}).
The second test is based on extensions of this and is known as the
Blum-Kiefer-Rosenblatt's (BKR) independence test; see~\cite{BKR}. 
To perform the BKR test, one can use their test statistic, $B_n$, together with the normal approximation to its null distribution as discussed in \cite{MW}.
Since the results of these two statistics
are always quite similar, here  we only report the results based on Hoeffding's statistic.

Another well known test for independence in the case $d=2$  is based on  Spearman's $\rho$, and  has been used extensively in applications.  However, it is well known that this test has low power if the distribution under the alternative hypothesis is continuous but  singular, as is the case for several copulas.
We used a small value of the sample size,  $n=36$, to compare our results
to those of other papers that also use small sample sizes in their simulations.
We made use of the {\bf spearman.test} function of the {\sc R package pspearman}.

Figure~\ref{Clayton} shows the power comparisons for the Clayton family. We observe that the power obtained by using the tests of Hoeffding, Blum-Kiefer-Rosenblatt and Spearman's $\rho$ are slightly better at levels $\alpha=0.01, 0.05$ and $\alpha=0.10$ than the ones we obtain using our statistics based on the total variation and Hellinger distances, and the Kullback-Leibler divergence. This behavior can also be observed for other standard Archimedean copulas such as Gumbel and Frank, among others. It is important to note that most of all these copulas are absolutely continuous, with complete support and with smooth densities.
We did not use the supremum distance in the simulations because we observed a strong discretization effect on the distribution of the statistic (\ref{statsup}) when the sample size is small. In other words, the different possible values of this statistic are very limited, leading to many ties in the simulated values.

It is not difficult to see, via simulations, that the independent tests of Hoeffding and Blum-Kiefer-Rosenblatt can also have problems with small sample sizes. In fact, we noted that when we sample from the independent copula $\Pi_2$, and test at the usual levels $\alpha \in \{0.01, 0.05, 0.10\}$, the real levels of the test do not correspond to the desired values of $\alpha$.
For example, if we set $\alpha=0.05$ and  perform several simulations, the
actual value of $\alpha$ under independence is approximately $0.075$.
Something similar happens with the other two values of $\alpha$.
This is due to the discrete nature of the test statistic, and the problem is more severe when the sample size is small. For this reason, we recommend caution when using these two tests with small sample sizes.

As a second example, we used the Fr\'echet-Mardia copulas. In this case, 
we use a convex mixture of $W_2$ and $M_2$, the Fr\'echet-Hoeffding bounds. As we can see in Figure~\ref{Fre-Mar}, the Spearman's test has very low power, specially for  singular copulas. We also note that the total variation statistic in equation (\ref{stattvd}) performs a little better than the Hoeffding and Blum-Kiefer-Rosenblatt tests at the three levels, but our test statistics based on the Hellinger distance and Kullback-Leibler divergence have the best performance at all three levels, and have a power close to 100\% when $\alpha=0.05$ and $\alpha=0.10$.

Finally, we used a convex combination of a Gumbel and a Gumbel-ID, where the latter denotes a Gumbel distribution with an increasing transformation in its first coordinate and a decreasing transformation in its second coordinate. That is, we applied the transformation $(U,V)\rightarrow (U,1-V)$ to all the observations $(U_i,V_i)$ from the Gumbel copula. The notation $\mbox{ID}$ stands for {\it increasing-decreasing} transformation. As we can see in Figure~\ref{Gum-GumID}, in this case the Spearman's $\rho$, the Hoeffding and the Blum-Kiefer-Rosenblatt have lower powers than our three test statistics (\ref{stattvd}), (\ref{statHel}) and (\ref{statKL}). Moreover, the performance of these three statistics can be much better than the ones obtained using the standard tests. In particular, the Kulback-Leibler test (\ref{statKL}) has the highest powers.

Summing up, our statistics (\ref{stattvd}), (\ref{statHel}) and (\ref{statKL}) are competitive in the case $d=2$. Also, for very smooth copulas with complete support and $0.5<|\hat{\rho}|< 0.95$ (where $\hat{\rho}$ is the estimated value of Spearman's $\rho$) we found that Spearman's test had the best powers for small sample sizes.

\end{subsection}

\begin{subsection}{Dimension $d=3$}

Many of the test statistics proposed in dimension $d=2$ can be extended to higher dimensions. Some of these extensions are rather natural. For example, in the case of the Hoeffding statistic and dimension $d=3$, we could use
$D(x,y,z)=F(x,y,z)-F_1(x)\cdot F_2(y)\cdot F_3(z)$, where
$F$ denotes the joint distribution function and $F_1$, $F_2$ and $F_3$
are the margins of $X$, $Y$ and $Z$, respectively. If we now define
$\Delta(F)= \int D^2(x,y,z) dF(x,y,z)$, we could use the empirical version of $\Delta$ to test for independence. The Blum-Kiefer-Rosenblatt statistic could  be similarly extended. Many other test statistics based on the empirical copula have
$3$-dimensional versions, for example the statistic $G_n$ of Genest and R\'emillard \cite{GR}.  The  problem with all these possible extensions is that in dimension $d=3$ the empirical copula may become difficult to compute if the sample size is not  small. 
On the other hand, the use of the empirical distribution function in dimensions greater than or equal to $d=3$ requires large sample sizes in order to obtain reasonable approximations of the true distribution function. In the case of $3$-copulas, the same is true of the empirical (sub)copula.

In some cases, one can find the asymptotic distribution of a test statistic based on the empirical distribution function or empirical copula, but this limiting  distribution can only be reached with large sample sizes. In such cases, the test statistic may be difficult to evaluate, and so it would not be possible to assess for which sample sizes the limiting distribution actually provides a good approximation.

There are other tests, based on empirical processes or multivariate characteristic functions, which also have problems when working with large sample sizes; see, for example, \cite{FLPS}. In \cite{BRS}, the authors propose a statistic, $\cali_{n}^{2}$, which coincides with the square product moment correlation when $d=2$. The power of this test is adequate only for absolutely continuous random variables, and it has the same problem as Spearman's test in dimension $d=2$. On the other hand, in \cite{TRO}, the authors propose a new test of multivariate independence based on analogues of Kendall's $\tau$ and Spearman's $\rho$. The comments of the previous paragraph also apply to these tests.

In the simulations using the four test statistics given by equations (\ref{stattvd}), (\ref{statHel}), (\ref{statsup}) and (\ref{statKL}), we used relatively large values of the sample size $n$. In many instances the other tests take a very long time.  Therefore, we only have compared our proposals among themselves in order to see which one has better powers in a number of scenarios.

In Figures~\ref{Gumbel} through~\ref{Normd30R} we analyze the case $d=3$. We used three sample sizes, $n=60$, $n=120$ and $n=216$, with $N=10\, 000$ simulations to find the critical values of the tests under $H_0$.  We also generated $1\, 000$ simulations  under the alternative hypothesis, $H_1$, to compute the powers of the four tests.

In Figure~\ref{Gumbel} we consider the Gumbel family, with $\alpha = 0.10$. We observe that, for small $n$, the test statistic based on the supremum distance has better powers, while for $n=216$ the powers are similar for all tests. As pointed out before, the statistic based on the supremum distance is highly discrete, which means that the critical values for this test are not very accurate. The same can be said about Figure~\ref{Normd3}, where we consider the normal family with a covariance matrix that has the same correlation for each pair of variables. We also note that the supremum distance has a little better power, but this advantage dissipates as the sample size increases. In Figure~\ref{Normd30R} we also study the normal family, but now one of the variables is independent of the other two. Note that in this case the test based on the supremum distance has the worst power compared to the other three statistics, but for large values of $n$ this difference seems to disappear.
\end{subsection}

\begin{subsection}{Dimension $d=4$}

The comments made for the case $d=3$ also apply to the case $d=4$. The only difference is that here we consider different values of the sample size which now include $n=600$ and $n=1296$. (The value $1296=(16)\cdot(81)$ is obtained by multiplying the number of boxes of $C^{(2)}$ and $C^{(3)}$ in dimension $d=4$.)
It is worth emphasizing that, if we tried to evaluate the empirical distribution function of a sample of size $n=1296$ in 4 dimensions, we would most likely get an error message because the array needed to store it would be of size $(1296)^4=2\, 821\, 109\, 907\, 456$, which no standard personal computer can handle. 

In Figure~\ref{Frank} we consider the Frank family with $\alpha = 0.05$. The remarks for this case are similar to those relating to Figure~\ref{Gumbel}. Finally, in Figure~\ref{St4} we study the Student $t$ distribution with $4$ degrees of freedom and having the same correlation for all the pairs of variables. As is well known, this distribution has heavier tails than the normal distribution, and in this case the test based on the Kullback-Leibler divergence performs better than the other tests.
\end{subsection}

\end{section}

\begin{figure}[h!]
\begin{center}
\includegraphics[width=5.1cm,height=5.1cm]{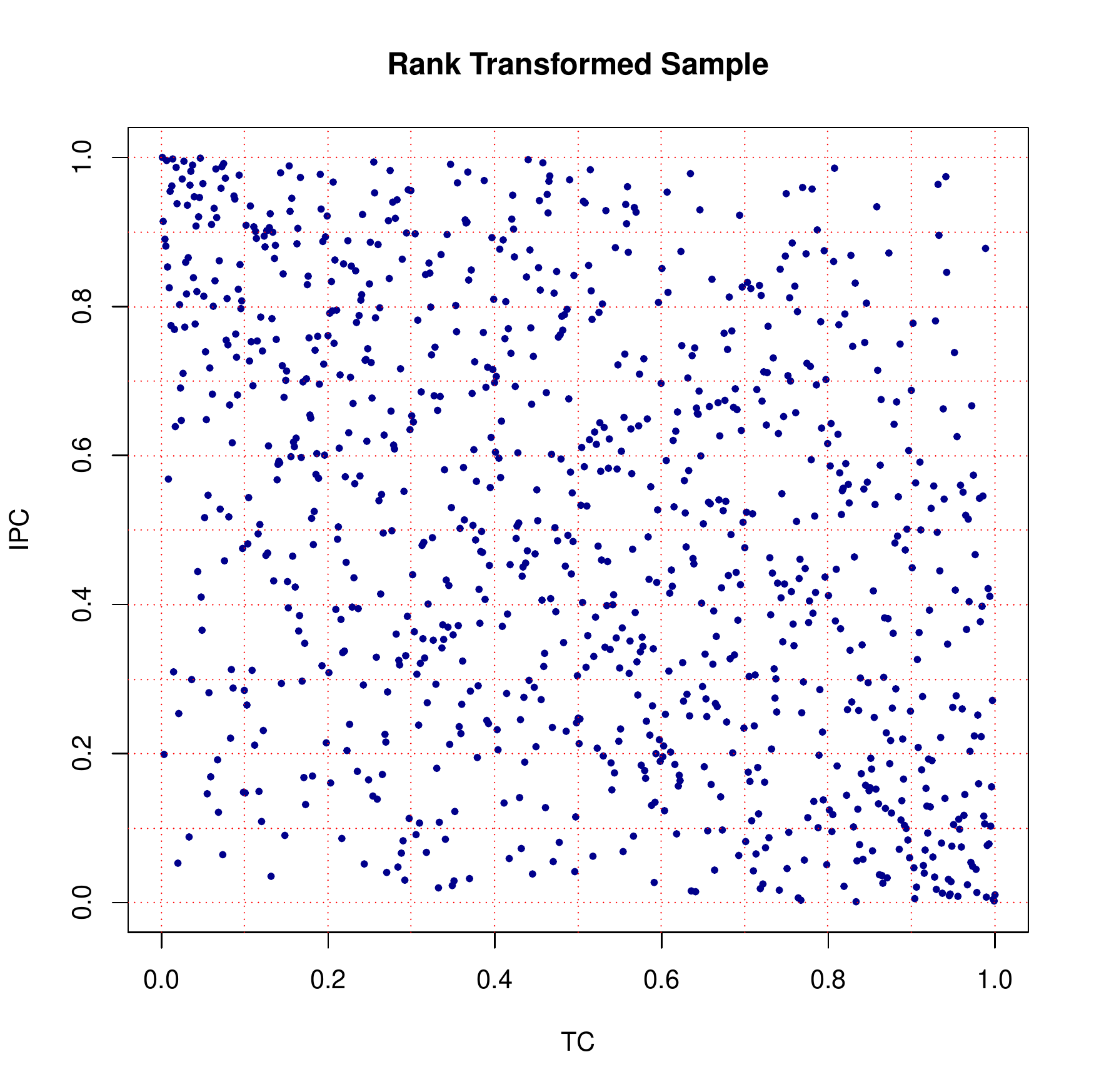}\hspace{.1cm}
\includegraphics[width=5.1cm,height=5.1cm]{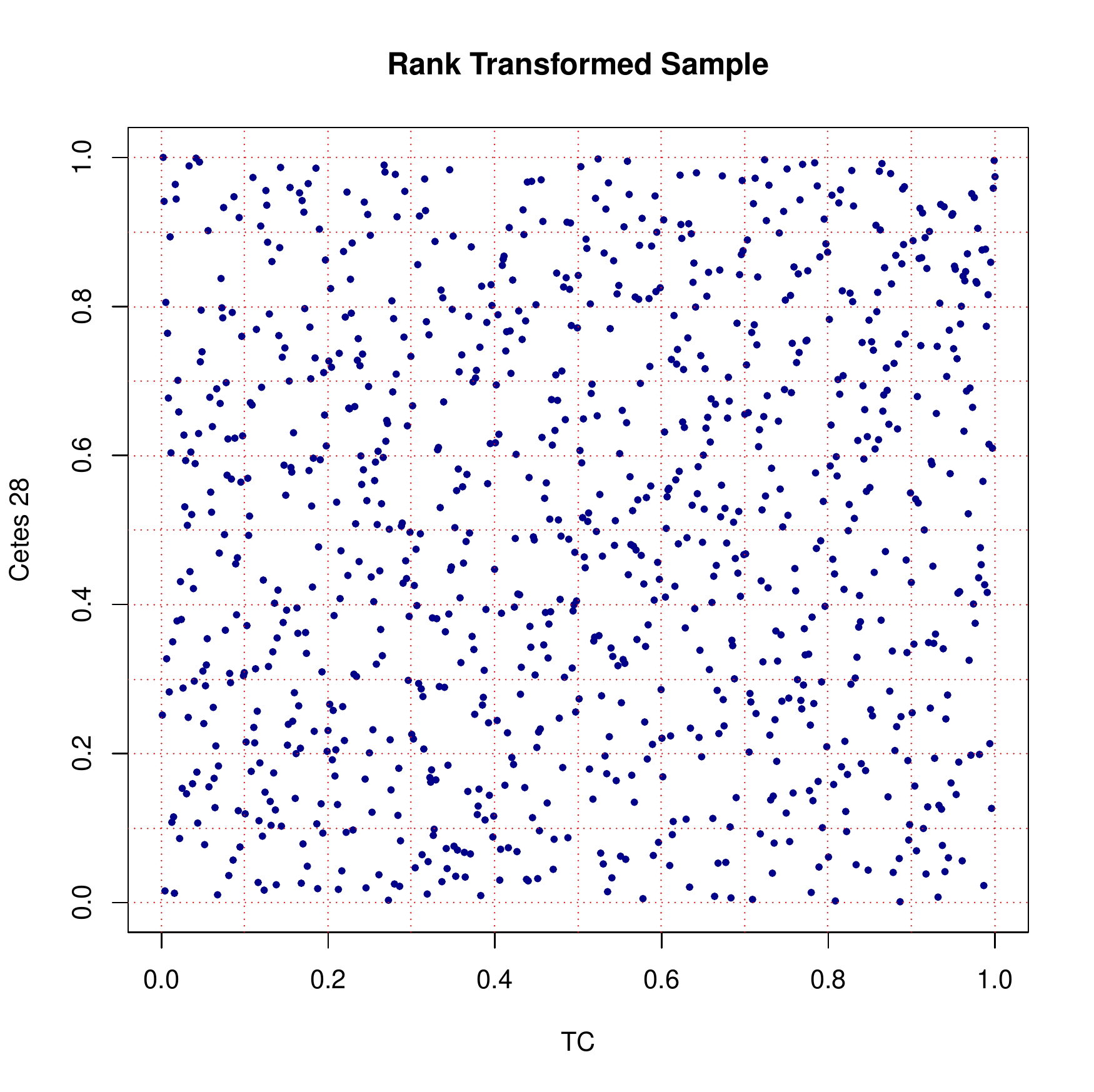}\hspace{.1cm}
\includegraphics[width=5.1cm,height=5.1cm]{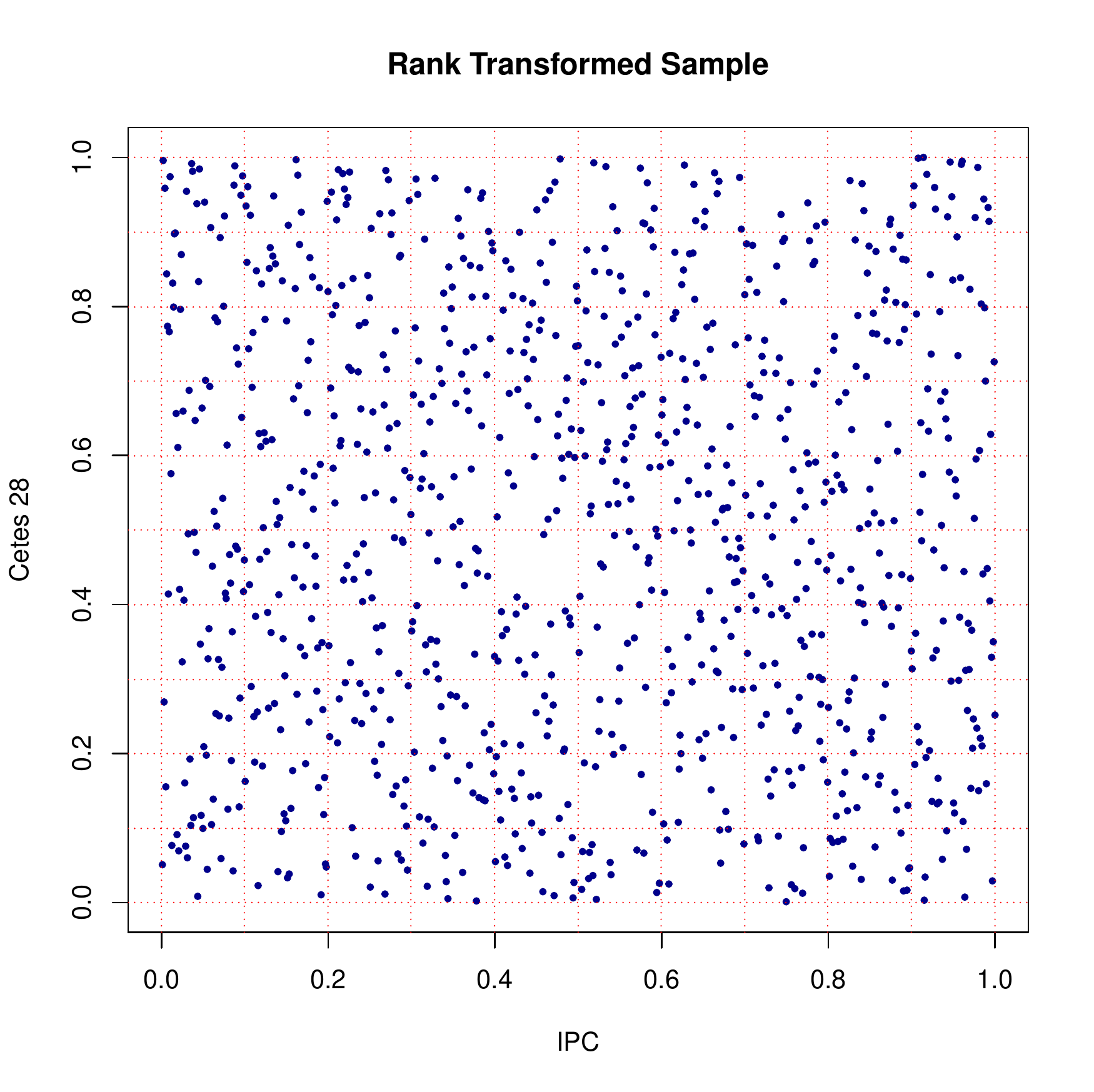}\hspace{.1cm}
\caption{TC-IPC-Cetes 28 real data set: pairwise scatter plots of the transformed samples.}
\label{RealData}
\end{center}
\end{figure}

\begin{section}{Real data example}

We applied our independence tests to real Mexican economic data in dimensions~$2$ and~$3$. We used data recorded on $n=967$ consecutive days, from 2014 to 2017, concerning three variables: the USD/MXN Exchange Rate (\textit{Tipo de Cambio} in Spanish, denoted here by TC); the Prices and Quotations Index of the Mexican Stock Exchange (\textit{Indice de Precios y Cotizaciones de la Bolsa Mexicana de Valores}, denoted here by IPC); and the price of a Mexican bond known as {\it Cetes 28}, where 28 refers to days. As is common in financial contexts, we did not worked with the raw data but used the corresponding returns instead.

Figure~\ref{RealData} shows the scatter plots of the modified (rank transformed) returns for the pairs (TC,IPC), (TC,Cetes 28) and (IPC,Cetes 28). The corresponding values of the Pearson correlation coefficients are $r = -0.3792$ for (TC,IPC), $r = 0.0944$ for (TC,Cetes 28) and $r = -0.0321$ for (IPC,Cetes 28).
When we applied our independence tests to the 3-dimensional data set, we rejected independence with all of them; that is, using the test statistics based on the total variation, supremum and Hellinger distances, as well as the Kullback-Leibler divergence.

We also applied our independence tests, together with Hoeffding, Blum-Kiefer-Rosenblatt and Spearman tests, to each of the three pairs of variables. For the first two pairs, (TC,IPC) and (TC,Cetes 28), all of the tests rejected independence at levels $0.10$, $0.05$ and $0.01$. However, for the pair (IPC, Cetes 28) none of the tests rejected independence at any of the three levels.

\end{section}

\begin{section}{Discussion}

In this paper we have provided a simple characterization of multivariate independence in terms of the checkerboard approximations of order~2 and~3 to a $d$-variate copula. While interesting in its own right, this result has also allowed us to propose a new family of tests of multivariate independence for $d$-variate continuous random vectors.

Our test statistics are all functionals of the sample copulas that estimate the above-mentioned checkerboard approximations.
These estimators can be evaluated for relatively large sample sizes even if the dimension is not small. This allows us to produce a large number of simulations in order to estimate the null distributions of the statistics we propose.
On the other hand, our test statistics can be defined in terms of any metric, or ever in terms of other functionals that are not symmetric such as the Kullback-Leibler divergence.

We simulated a range of examples in dimensions up to $d=4$, under different
models and with different sample sizes.
In many of these scenarios, it may not be feasible to compute the empirical distribution functions. 
In our simulations we observed that, when the sample size is moderately
large, all of the tests we considered have similar powers. Thus, in this case 
any of our test statistics may be used. However, when the sample size is small, we warn the user against the test statistic based on the supremum distance since it is strongly affected by discretization.

One interesting application of our tests is the following. Consider a random vector ${\bf X}=(X_1,X_2,\ldots ,X_d)$. In some cases it is possible to decompose ${\bf X}$ into two independent subvectors ${\bf Y_1}$ and ${\bf Y_2}$; that is, we can find a permutation $\pi$ of
$\{1,2,\ldots ,d-1\}$ and a value of $k\in \{1,2,\ldots,d\}$ such that the subvectors
${\bf Y_1}=(X_{\pi(1)},X_{\pi(2)},\ldots, X_{\pi(k)})$ and ${\bf Y_2}=(X_{\pi(k+1)},\ldots ,X_{\pi(d)})$ are independent. Conversely, suppose there exist no permutation $\pi$ of
$\{1,2, \ldots ,d\}$ and $k\in\{1,2,\ldots ,d\}$ such that, for every 
$(x_{\pi(1)},\ldots ,x_{\pi(d)})\in\rea^d$,
\begin{equation}
	F_{\pi(1),\ldots ,\pi(d)}(x_{\pi(1)},\ldots ,x_{\pi(d)})=
	F_{\pi(1),\ldots ,\pi(k)}(x_{\pi(1)},\ldots ,x_{\pi(k)})\cdot 
	F_{\pi(k+1),\ldots ,\pi(d)}(x_{\pi(k+1)},\ldots ,x_{\pi(d)})
\end{equation}
where $F_{\pi(1),\ldots ,\pi(k)}$ and $F_{\pi(k+1),\ldots ,\pi(d)}$ are marginal distribution functions. In this case we follow \cite{GBGPR} and say that the random vector ${\bf X}$ is \textbf{exhaustively dependent}.

In practice, it is not uncommon to find random vectors for which it can be assumed that a certain set of coordinates are independent from the remaining coordinates. In such cases, for  a  sample of size $n$ of this random vector, we would like to produce a statistical test for the independence of these two subvectors. 
In \cite{GBGPR} the authors show that, if 
${\bf X}=(X_1,X_2,\ldots ,X_d)$ is a random vector on $\rea^d$ ($d\geq 3$) with joint distribution function
$F_{1,2,\ldots ,d}$, and if there exists a permutation $\pi$ of $\{1,2,\ldots,d\}$ such that all the subvectors
$(X_{\pi(1)},X_{\pi(2)})$, $(X_{\pi(2)},X_{\pi(3)}), \ldots,(X_{\pi(d-1)},X_{\pi(d)})$ are dependent, then
${\bf X}=(X_1,X_2,\ldots ,X_d)$ is exhaustively dependent.

Assume that $d=5$  and take $\pi$ to be the identity. Now let ${\bf X}=(X_1,X_2,X_3,X_4,X_5)$ be a random vector such that
${\bf Y_1}=(X_1,X_2,X_3)$ and ${\bf Y_2}=(X_4,X_5)$ are two subvectors which are exhaustively dependent, and assume that we suspect that ${\bf Y_1}$ and ${\bf Y_2}$ are independent. If we have a random sample of size $n$
from the distribution of ${\bf X}$, we can verify all these assumptions
using our tests of independence as follows.
First, test if the subvectors $(X_1,X_2) , (X_2,X_3)$ for ${\bf Y_1}$
and $(X_4,X_5)$ for ${\bf Y_2}$
have independent coordinates (at a certain level $\alpha$ for each test). If we reject the hypothesis of independence in all three tests,
now we can test independence of the subvectors ${\bf Y_1}$ and ${\bf Y_2}$ by testing the independence between the coordinates of each of the 
following six subvectors: $(X_1,X_4), (X_2,X_4), (X_1,X_5), (X_2,X_5), (X_3,X_4)$ and $(X_3,X_5)$. If we do not reject the hypothesis of 
independence for all these pairs, then we have some evidence that our assumptions are not incorrect. Note that in this case all the independence tests are performed for pairs of random variables, which makes the tests quite quick even for large sample sizes.

We have developed a program in the statistical language R that implements the procedure described above. The code is available from the authors upon request.

\end{section}

\section*{Appendix}

\subsection*{Proof of Theorem~\ref{MainThm}}

First, assume that $d=2$. Let us assume that $C=\Pi_2$, that is, $C$ is the independence copula; then, using equation (\ref{cm}), we have

$$ C_{m}^{'} =\{ C(u,v)=u\cdot v \,\,|\,\, u,v\in\{ 0,1/m,2/m, \ldots ,(m-1)/m,1\} \},$$
is a 2-subcopula. For this 2-subcopula and the uniform partition of size $m$ given in equation
(\ref{partunif}), and using equation (\ref{volumen}), we have that for every $i_1, i_2 \in I_m$,

\begin{eqnarray}\label{lebesgue}
V_{C_{m}^{'}}(\overline{R_{i_1,i_2}^{m}}) &=&
\frac{i_1}{m}\frac{i_2}{m}-\frac{i_1-1}{m}\frac{i_2}{m}-\frac{i_1}{m}\frac{i_2-1}{m}+
\frac{i_1-1}{m}\frac{i_2-1}{m}\nonumber\\
& = & \left(\frac{i_1}{m}-\frac{i_1-1}{m}\right) \left(\frac{i_2}{m}-\frac{i_2-1}{m}\right)\nonumber\\
& = & \lambda^2(\overline{R_{i_1,i_2}^{m}}),
\end{eqnarray}
where $\lambda^2$ is the Lebesgue measure on $(\rea^2, \calb(\rea^2))$.

If we use the bilinear interpolation  of Lemma 2.3.5 in Nelsen's book, see \cite{N06},
we have that $C^{(m)}$ the checkerboard approximation of order $m$ of $C=\Pi_2$ has 
a density given by equation (\ref{denschecker})

\begin{equation}\label{density}
c^{(m)}(u,v) =\frac{V_{C_{m}^{'}}
(\overline{R_{i_1,i_2}^{m}})}{\lambda^2(\overline{R_{i_1,i_2}^{m}})}
\quad\mbox{for every}\quad (u,v)\in R_{i_1,i_2}^{m},
\end{equation}
for every $i_1,i_2\in I_m$. On the other hand, using equations (\ref{lebesgue}) and (\ref{density}) we have that

$$ c^{(m)}(u,v) =\frac{V_{C_{m}^{'}}(\overline{R_{i_1,i_2}^{m}})}
{\lambda^2(\overline{R_{i_1,i_2}^{m}})}=
\frac{\lambda^2(\overline{R_{i_1,i_2}^{m}})}{\lambda^2(\overline{R_{i_1,i_2}^{m}})}
= 1 \,\,\mbox{for every}\,\, (u,v)\in R_{i_1,i_2}^{m},$$
for every $i_1,i_2\in I_m$. Hence, the density of $C^{(m)}$ is the constant $1$ on ${\bf I}^2$.
Therefore, for every integer $m\geq 2$, the checkerboard approximation $C^{(m)}$ satisfies 
\begin{equation}\label{gral}
C^{(m)}(u,v)=\int_{0}^{v}\int_{0}^{u} 1 ds dt=u\cdot v=\Pi_{2}(u,v)=C(u,v)\,\,\mbox{for every}\,\,(u,v)\in {\bf I}^2.
\end{equation}
In particular this holds for $m=2$ and $m=3$.

Now, let us assume that for some 2-copula $C$ we have that
$C(u,v)=C^{(2)}(u,v)=C^{(3)}(u,v)$ for every $(u,v)\in {\bf I}^2$.
Let $m=2$ and define $\alpha=V_C([0,1/2]^2)=V_C(R_{1,1}^2)$, as  in the uniform partition of order $m=2$, given in
equation (\ref{partunif}). Then, by equation (\ref{volumen}) and using inequality (\ref{FreHoeff}), if $\alpha =C(1/2,1/2)$, we have
\begin{equation}\label{bounds}
0=W(1/2,1/2)\leq \alpha=C(1/2,1/2)\leq M(1/2,1/2)=\frac{1}{2}.
\end{equation}
Note that $R_{1,1}^{2}\cup R_{1,2}^{2} =[0,1/2]\times[0,1]$ is a disjoint union. 
Note also that, by continuity of $C$, $V_C(R_{1,2}^{2})=V_C(\overline{R_{1,2}^{2}})$; the same applies to
$R_{2,1}^{2}$ and $R_{2,2}^{2}$. Hence, using equation (\ref{volumen}),
$$
\frac{1}{2}=V_C([0,1/2]\times[0,1])=V_C(R_{1,1}^{2})+V_C(R_{1,2}^{2})=\alpha +V_C(R_{1,2}^{2})
$$
and so $V_C(R_{1,2}^{2})=(1/2-\alpha)$. Similar arguments show that $V_C(R_{2,1}^{2})=(1/2-\alpha)$ and that
$V_C(R_{2,2}^{2})=\alpha$. So, using the bilinear interpolation we obtain
\begin{equation}\label{c2}
C^{(2)}(u,v)=C_{\alpha}(u,v)=\left\{\begin{array}{lcl}
4\alpha uv  & \mbox{if} & (u,v)\in R_{1,1}^{2} \nonumber\\
2\alpha u +4(1/2-\alpha)u(v-1/2) & \mbox{if} & (u,v)\in R_{1,2}^{2} \nonumber\\
2\alpha v +4(1/2-\alpha)(u-1/2)v & \mbox{if} & (u,v)\in R_{2,1}^{2} \nonumber\\
\alpha +(1-2\alpha)(u+v-1)+4\alpha(u-1/2)(v-1/2) & \mbox{if} & (u,v)\in R_{2,2}^{2}.
\end{array}\right.
\end{equation}
From equation (\ref{c2}) we have that $C^{(2)}$ is a function of a unique parameter, that is,
$\alpha=C(1/2,1/2)$, and from the hypothesis we have that $C(u,v)=C^{(2)}(u,v)=C_{\alpha}(u,v)$,
where from equation (\ref{bounds}), $0\leq \alpha\leq 1/2$.

Now, we also assume that $C$ satisfies

\begin{equation}\label{c3} 
C_{\alpha}(u,v)=C(u,v)=C^{(2)}(u,v)=C^{(3)}(u,v)=C_{\alpha}^{(3)}(u,v),
\end{equation} 
for every $(u,v)\in{\bf I}^2$. In order to construct $C_{\alpha}^{(3)}(u,v)$ we need to evaluate
all the volumes $V_{C_{\alpha}}(\overline{R_{i_1,i_2}^{3}})$ for every $i_1,i_2\in I_3=\{1,2,3\}$.
We first observe that $R_{1,1}^3=[0,1/3]^2\subset [0,1/2]^2=R_{1,1}^2$, so using equation (\ref{c2}),
we obtain
\begin{equation}\label{volr11}
V_{C_{\alpha}}(R_{1,1}^{3})=V_{C_{\alpha}}([0,1/3]^2)=C_{\alpha}(1/3,1/3)=\frac{4\alpha}{9}.
\end{equation}
In general, by continuity of $C$, $V_{C_{\alpha}}(R_{i_1,i_2}^{3})
=C_{\alpha}(i/m,j/m)-C_{\alpha}((i-1)/m,j/m)-C_{\alpha}(i/m,(j-1)/m)+C_{\alpha}((i-1)/m,(j-1)/m)$
for every $i_1,i_2\in I_3$. We also know from equation (\ref{partunif}), that 
$\lambda^2(R_{i_1,i_2}^{3})=1/9$ for every $i_1,i_2\in I_3$.
Hence, using equation (\ref{density}), the density of $C_{\alpha}^{(3)}$ is given by
\begin{equation} \label{densityc3}
c_{\alpha}^{(3)}(u,v)=\frac{V_{C_{\alpha}}(\overline{R_{i_1,i_2}^{3}})}
{\lambda^2(\overline{R_{i_1,i_2}^{3}})}= 9V_{C_{\alpha}}(\overline{R_{i_1,i_2}^{3}})
\end{equation}
for every $(u,v)\in R_{i_1,i_2}^{3}$ and for every $i_1,i_2\in I_3$. Using equations (\ref{c2}) and (\ref{volr11}) we have that

\begin{equation}\label{valuer11}
C_{\alpha}^{(3)}(u,v) = 9V_{C_{\alpha}}(R_{1,1}^{3})u\cdot v = 
9\left(\frac{4\alpha}{9}\right)u\cdot v = 4\alpha u\cdot v,
\end{equation}
for every $(u,v)\in R_{1,1}^{3}=[0,1/3]^2$. We also have that
\begin{eqnarray}\label{volr12}
V_{C_{\alpha}}(R_{1,2}^{3}) &=& C_{\alpha}(1/3,2/3)-C_{\alpha}(1/3,1/3)-C_{\alpha}(0,2/3)+C_{\alpha}(0,1/3)\nonumber\\
& = & \frac{2\alpha}{3} +4\left(\frac{1}{2}-\alpha\right) \left(\frac{1}{3}\right)\left(\frac{1}{6}\right) -\frac{4\alpha}{9} \nonumber\\
& = & \frac{1}{9} +\alpha \cdot\frac{12-4-8}{18}\nonumber\\
& = & \frac{1}{9}
\end{eqnarray}
Now, using equations(\ref{valuer11}), (\ref{volr12}) and integration we obtain
\begin{equation}\label{valuer12}
C_{\alpha}(u,v)=9V_{C_{\alpha}}(R_{1,1}^{3})u\left(\frac{1}{3}\right) +9V_{C_{\alpha}}(R_{1,2}^{3})u\left(v-\frac{1}{3}\right),
\end{equation}
for every $(u,v)\in R_{1,2}^{3}=[0,1/3]\times(1/3,2/3]$.

Finally, let us take $(u_0,v_0)=(1/4,1/2)$. Using equation (\ref{partunif}), we have
that $(1/4,1/2)\in (R_{1,1}^{2}\cap R_{1,2}^3)$, while from equation (\ref{c2}) we have
\begin{equation}\label{c2u0v0}
C_{\alpha}(u_0,v_0)=C^{(2)}(1/4,1/2)=4\alpha \left(\frac{1}{4}\right) \left(\frac{1}{2}\right)=\frac{\alpha}{2}.
\end{equation}
On the other hand, from equations (\ref{volr12}) and (\ref{valuer12}) it follows that

\begin{equation}\label{c3u0v0}
C_{\alpha}^{(3)}(1/4,1/2)=4\alpha \left(\frac{1}{4}\right) \left(\frac{1}{3}\right) 
+9 \left(\frac{1}{9} \right) \left(\frac{1}{4}\right) \left(\frac{1}{6}\right) =\frac{\alpha}{3} +\frac{1}{24}.
\end{equation}
Therefore, from hypothesis (\ref{c3}) and equations (\ref{c2u0v0}) and (\ref{c3u0v0}), we have that
$$ C_{\alpha}^{(3)}(1/4,1/2)=C_{\alpha}(1/4,1/2)\,\,\, \mbox{if and only if}\,\,\,
\alpha/2=\alpha/3+1/24 \,\,\,\mbox{if and only if}\,\,\, \alpha=1/4.$$  
However, by equation (\ref{c2}), this happens if and only if $C_{\alpha}(u,v)=\Pi_2(u,v)=u\cdot v$. 

We now assume that $d=3$ and that $C=\Pi_3$ is the product 3-copula.
Then using equation (\ref{cm}) we know that
$$
C_{m}^{'} =\{ C(u,v,w)=u\cdot v \cdot w\,\,|\,\, u,v,w\in\{ 0,1/m,2/m, \ldots ,(m-1)/m,1\} \},
$$
is a 3-subcopula, and for this 3-subcopula and the uniform partition of size $m$ given in equation
(\ref{partunif}), and using equation (\ref{volumen}), we have by continuity of $C$ that
\begin{eqnarray}\label{lebesgue3}
V_{C_{m}^{'}}(R_{i_1,i_2,i_3}^{m}) &=&
\frac{i_1}{m}\frac{i_2}{m}\frac{i_3}{m}-\frac{i_1-1}{m}\frac{i_2}{m}\frac{i_3}{m}
-\frac{i_1}{m}\frac{i_2-1}{m}\frac{i_3}{m}-\frac{i_1}{m}\frac{i_2}{m}\frac{i_3-1}{m} \nonumber\\
&  & +\frac{i_1-1}{m}\frac{i_2-1}{m}\frac{i_3}{m} +\frac{i_1-1}{m}\frac{i_2}{m}\frac{i_3-1}{m}
+\frac{i_1}{m}\frac{i_2-1}{m}\frac{i_3-1}{m}-\frac{i_1-1}{m}\frac{i_2-1}{m}\frac{i_3-1}{m}\nonumber\\
& = & \left(\frac{i_1}{m}-\frac{i_1-1}{m}\right) 
      \left(\frac{i_2}{m}-\frac{i_2-1}{m}\right)\left(\frac{i_3}{m}-\frac{i_3-1}{m}\right) \nonumber\\
& = & \lambda^3(R_{i_1,i_2,i_3}^{m}),
\end{eqnarray}
where $\lambda^3$ is the Lebesgue measure on $(\rea^3, \calb(\rea^3))$.

If we use the trilinear interpolation  of Lemma 2.3.5 in Nelsen's book, \cite{N06},
we have that $C^{(m)}$ the checkerboard approximation of order $m$ of $C=\Pi_3$ has 
a density given by equation (\ref{denschecker})
\begin{equation}\label{density3}
c^{(m)}(u,v,w) =\frac{V_{C_{m}^{'}}(R_{i_1,i_2,i_3}^{m})}{\lambda^3(R_{i_1,i_2,i_3}^{m})}
\quad\mbox{for every}\quad (u,v,w)\in R_{i_1,i_2,i_3}^{m},
\end{equation}
for every $i_1,i_2,i_3\in I_m$. But, using equations (\ref{lebesgue3}) and (\ref{density3}) we have that

$$ c^{(m)}(u,v,w) =
\frac{V_{C_{m}^{'}}(R_{i_1,i_2,i_3}^{m})}{\lambda^3(R_{i_1,i_2,i_3}^{m})}=
\frac{\lambda^3(R_{i_1,i_2,i_3}^{m})}{\lambda^3(R_{i_1,i_2,i_3}^{m})}
= 1 \,\,\mbox{for every}\,\, (u,v,w)\ R_{i_1,i_2,i_3}^{m},$$
for every $i_1,i_2,i_3\in I_m$. Hence, the density of $C^{(m)}$ is the constant $1$ on ${\bf I}^3$.
Therefore, for every integer $m\geq 2$ the checkerboard approximation $C^{(m)}$ satisfies that
\begin{equation}\label{gral3}
C^{(m)}(u,v,w)=\int_{0}^{v}\int_{0}^{u}\int_{0}^{w} 1 \, ds dt dr
=u\cdot v\cdot w=\Pi_{3}(u,v,w)=C(u,v,w)\,\,\mbox{for every}\,\,(u,v,w)\in {\bf I}^3.
\end{equation}
In particular this holds for $m=2$ and $m=3$.

We now prove the converse. Let us assume that for some 3-copula $C$ we have that
$C(u,v,w)=C^{(2)}(u,v,w)=C^{(3)}(u,v,w)$ for every $(u,v,w)\in {\bf I}^3$.
Let $m=2$ and define $\alpha_0=V_C([0,1/2]^3)=V_C(R_{1,1,1}^2)$ as in the uniform partition of order $m=2$ given in
equation (\ref{partunif}). Then, by equation (\ref{volumen}), $\alpha_0 =C(1/2,1/2,1/2)$;
using now the inequality (\ref{FreHoeff}), we have
\begin{equation}\label{bounds3}
0=W^3(1/2,1/2,1/2)\leq \alpha_0=C(1/2,1/2,1/2)\leq M^3(1/2,1/2,1/2)=\frac{1}{2}.
\end{equation}
Define $\alpha_1=C(1,1/2,1/2), \alpha_2=C(1/2,1,1/2)$ and $\alpha_3=C(1/2,1/2,1)$.
Let $C_{1,2}(u,v)=C(u,v,1)$, the we know that $C_{1,2}$ is a 2-copula, and by hypotheses
we also know that $C_{1,2}(u,v)=C^{(2)}(u,v,1)=C^{(3)}(u,v,1)$ for every $(u,v)\in {\bf I}^2$.
It is trivial to see that by linearity in the construction of $C^{(2)}$ and $C^{(3)}$, we have that
the checkerboards of $C_{1,2}$ of order $m=2$ and $m=3$ are given by
$C_{1,2}^{(2)}(u,v)=C^{(2)}(u,v,1)$ and $C_{1,2}^{(3)}(u,v)=C^{(3)}(u,v,1)$ for every $(u,v)\in{\bf I}^2$.
Therefore, we have the transformed hypotheses
\begin{equation}\label{hypinduction}
C_{1,2}(u,v)=C_{1,2}^{(2)}(u,v)=C_{1,2}^{(3)}(u,v)\quad\mbox{for every}\quad (u,v)\in{\bf I}^2.
\end{equation}
So, using what we proved in the case $d=2$ above,
\begin{equation}\label{alpha3}
\alpha_3=C(1/2,1/2,1)=C_{1,2}(1/2,1/2)=\Pi_2(1/2,1/2)=\frac{1}{4}.
\end{equation}
Defining $C_{1,3}(u,w)=C(u,1,w)$ and $C_{2,3}(v,w)=C(1,v,w)$ for every $u,v,w\in{\bf I}$, and
reasoning as above we observe that
\begin{equation}\label{alphas}
\alpha_1=C(1,1/2,1/2)=C_{2,3}(1/2,1/2)=\frac{1}{4}=C_{1,3}(1/2,1/2)=C(1/2,1,1/2)=\alpha_2.
\end{equation}
Now using the fact that any 3-copula is increasing in each coordinate, together with equations (\ref{alpha3}) and (\ref{alphas})
and inequality (\ref{bounds3}) we have
\begin{equation}\label{boundalfa0}
0\leq \alpha_0 =C(1/2,1/2,1/2)\leq \min(\alpha_1,\alpha_2,\alpha_3)=\frac{1}{4}.
\end{equation}

In order to find $C^{(2)}(u,v,w)$, we need first to evaluate the $C$-volumes of all the uniform
boxes $R_{i,j,k}^{2}$ for every $i,j,k\in I_2$, in order to find its density in each box, which is given
by the constant $V_C(\overline{R_{i,j,k}^2})/\lambda^3(R_{i,j,k}^2)=8\cdot V_C(\overline{R_{i,j,k}^2})$ for every
$i,j,k \in I_2=\{1,2\}$.
We know that $V_C(\overline{R_{1,1,1}^2})=C(1/2,1/2,/12)=\alpha_0$. By equation (\ref{volumen})
and using i) in Definition \ref{deftrivial}, we have that
$V_C(\overline{R_{2,1,1}^2})= V_C([1/2,1]\times[0,1/2]\times[0,1/2])=C(1,1/2,1/2)-C(1/2,1/2,1/2)=
\alpha_1-\alpha_0=1/4-\alpha_0$, similarly $V_C(\overline{R_{1,2,1}^2})=V_C(\overline{R_{1,1,2}^2})=1/4-\alpha_0$.
Again, by equation (\ref{volumen}) and using i) and ii) in Definition \ref{deftrivial}, we obtain
$V_C(\overline{R_{2,2,1}^2})=V_C([1/2,/1]\times[1/2,1]\times[0,1/2])=
C(1,1,1/2)-C(1,1/2,1/2)-C(1/2,1,1/2)+C(1/2,1/2,1/2)=1/2-\alpha_1-\alpha_2+\alpha_0=
1/2-1/4-1/4+\alpha_0=\alpha_0$, analogously, $V_C(\overline{R_{2,1,2}^2})=V_C(\overline{R_{1,2,2}^2})=\alpha_0$.
Finally, using  Definition \ref{deftrivial} we have that
$V_C(\overline{R_{2,2,2}^2})=V_C([1/2,/1]\times[1/2,1]\times[1/2],1)=
1-1/2-1/2-1/2+1/4+1/4+1/4-\alpha_0=1/4-\alpha_0$.

Therefore, integrating the above density we get 
$C^{(2)}(u,v,w)$ the checkerboard copula of order $m=2$, for every $(u,v,w)\in {\bf I}^3$,
which is given by:

\begin{equation}\label{c2-d3}
C^{(2)}(u,v,w)  =  \left\{\begin{array}{lcl}
8\alpha_0 u\cdot v\cdot w &\,\,\mbox{if}\,\, &(u,v,w)\in R_{1,1,1}^2 \nonumber\\
(2-8\alpha_0)u\cdot v\cdot w +(8\alpha_0-1) u\cdot v &\,\,\mbox{if}\,\, &(u,v,w)\in R_{1,1,2}^2 \nonumber\\
(2-8\alpha_0)u\cdot v\cdot w +(8\alpha_0-1) u\cdot w &\,\,\mbox{if}\,\, &(u,v,w)\in R_{1,2,1}^2 \nonumber\\
(2-8\alpha_0)u\cdot v\cdot w +(8\alpha_0-1) v\cdot w &\,\,\mbox{if}\,\, &(u,v,w)\in R_{2,1,1}^2 \nonumber\\
8\alpha_0 u\cdot v\cdot w +(1-8\alpha_0) u\cdot v    &                  &                       \nonumber\\
+(1-8\alpha_0)u\cdot w +(8\alpha_0-1) u                &\,\,\mbox{if}\,\, &(u,v,w)\in R_{1,2,2}^2 \nonumber\\
8\alpha_0 u\cdot v\cdot w +(1-8\alpha_0) u\cdot v    &                  &                       \nonumber\\
+(1-8\alpha_0)v\cdot w +(8\alpha_0-1) v                &\,\,\mbox{if}\,\, &(u,v,w)\in R_{2,1,2}^2 \nonumber\\
8\alpha_0 u\cdot v\cdot w +(1-8\alpha_0) u\cdot w    &                  &                       \nonumber\\
+(1-8\alpha_0)v\cdot w +(8\alpha_0-1) w                &\,\,\mbox{if}\,\, &(u,v,w)\in R_{2,2,1}^2 \nonumber\\
(1/2-2\alpha_0)\{(u-1/2)+(v-1/2)+(w-1/2)\}           &                  &                       \nonumber\\
4\alpha_0\{(u-1/2)(v-1/2)+(u-1/2)(w-1/2)\}           &                  &                       \nonumber\\
4\alpha_0 (v-1/2)(w-1/2) +\alpha_0                   &                  &                       \nonumber\\
(2-8\alpha_0)(u-1/2)(v-1/2)(w-1/2)                   &\,\,\mbox{if}\,\, &(u,v,w)\in R_{2,2,2}^2.
\end{array}\right.
\end{equation}
Note that by hypothesis $C^{(2)}(u,v,w)=C(u,v,w)$, and that
by equation (\ref{c2-d3}) it  has a unique parameter  $\alpha_0$. 

In order to obtain $C^{(3)}$, we will obtain its density using equation (\ref{c2-d3}), that is,
\begin{equation}\label{dens3d3}
c^{(3)}(u,v,w) =\frac{V_C(\overline{R_{i,j,k}^3})}{\lambda^3(R_{i,j,k}^3)}= 27 V_{C^{(2)}}(\overline{R_{i,j,k}^3}),
\end{equation}
for every $i,j,k\in I_3$ and for every $(u,v,w)\in R_{i,j,k}^3$, as defined in equation (\ref{partunif}).
To find the density of $C^{(3)}$ on $R_{1,1,1}^3=[0,1/3]^3$ we observe that $R_{1,1,1}^3\subset R_{1,1,1}^2$,
so, using (\ref{c2-d3}), $V_{C^{(2)}}(R_{1,1,1}^3)=C^{(2)}(1/3,1/3,1/3)= (8/27)\alpha_0$.
To obtain the density of $C^{(3)}$ on $R_{1,2,1}^3=[0,1/3]\times(1/3,2/3]\times[0,1/3]\subset R_{1,1,1}^2 \cup R_{1,2,1}^2$, we need 
$V_{C^{(2)}}(\overline{R_{1,2,1}^3})=C^{(2)}(1/3,2/3,1/3)-C^{(2)}(1/3,1/3,1/3)=(2-8\alpha_0)(2/27)+(8\alpha_0-1)(1/9)-8\alpha_0(1/27)=1/27$.
For the density of $C^{(3)}$ on $R_{1,1,2}^3=[0,1/3]\times[0,1/3]\times(1/3,2/3]\subset R_{1,1,1}^2 \cup R_{1,1,2}^2$ we need
$V_{C^{(2)}}(\overline{R_{1,1,2}^3})=C^{(2)}(1/3,1/3,2/3)-C^{(2)}(1/3,1/3,1/3)=(2-8\alpha_0)(2/27)+(8\alpha_0-1)(1/9)-8\alpha_0(1/27)=1/27$.
Finally, for the density of $C^{(3)}$ on $R_{1,2,2}^3=[0,1/3]\times(1/3,2/3]\times(1/3,2/3]\subset 
R_{1,1,1}^2 \cup R_{1,2,1}^2\cup R_{1,1,2}^2\cup R_{1,2,2}^2$ we need
$V_{C^{(2)}}(\overline{R_{1,2,2}^3})=C^{(2)}(1/3,2/3,2/3)-C^{(2)}(1/3,2/3,1/3)-C^{(2)}(1/3,1/3,2/3)+C^{(2)}(1/3,1/3,1/3)=
8\alpha_0(4/27)+(1-8\alpha_0)(4/9)+(8\alpha_0-1)(1/3)-(2-8\alpha_0)(4/27)-(8\alpha_0-1)(2/9)+8\alpha_0(1/27)=1/27$.
Hence, from equation (\ref{dens3d3}), we have that $C^{(3)}$ has density $1$ on $R_{1,1,2}^3, R_{1,2,1}^3$ and $R_{1,2,2}^3$, and density
$8\alpha_0$ on $R_{1,1,1}^3$.

Now let $(u_0,v_0,w_0)=(1/4,1/2,2)\in R_{1,1,1}^2\cap R_{1,2,2}^3$, Then, by hypothesis, $C^{(2)}(1/4,1/2,1/2)=C^{(3)}(1/4,1/2,1/2)$.
Integrating the density of $C^{(3)}$ we have
\begin{eqnarray}
C^{(3)}(1/4,1/2,1/2) & = & \int_{0}^{1/3} \int_{0}^{1/3}\int_{0}^{1/4}8\alpha_0 du dv dw + \int_{0}^{1/3} \int_{1/3}^{1/2}\int_{0}^{1/4} du dv dw \nonumber\\
                     &   & +\int_{1/2}^{1/3} \int_{0}^{1/3}\int_{0}^{1/4} du dv dw + \int_{1/3}^{1/2} \int_{1/3}^{1/2}\int_{0}^{1/4} du dv dw \nonumber\\
										 & = & (2/9)\alpha_0 +(1/72)+(1/72)+(1/144)\nonumber\\
										 & = & (2/9)\alpha_0 + (5/144),
\end{eqnarray}
and using equation (\ref{c2-d3}) we know that $C^{(2)}(1/4,1/2,1/2)=\alpha_0/2$. Therefore,
$$ \frac{\alpha_0}{2} =\frac{2}{9}\alpha_0 +\frac{5}{144}.$$
Solving for $\alpha_0$ we have that $\alpha_0=1/8$, and using equation (\ref{c2-d3}), we have that $C^{(3)}(u,v,w)=C^{(2)}(u,v,w)
=\Pi_2(u,v,w)$ for every $(u,v,w)\in{\bf I}^3$. The rest of the proof follows by induction. \qed


\newpage
\vspace*{-1ex}

\vspace*{-1ex}

\begin{figure}[t!]
\begin{center}
\includegraphics[width=5.1cm,height=5.1cm]{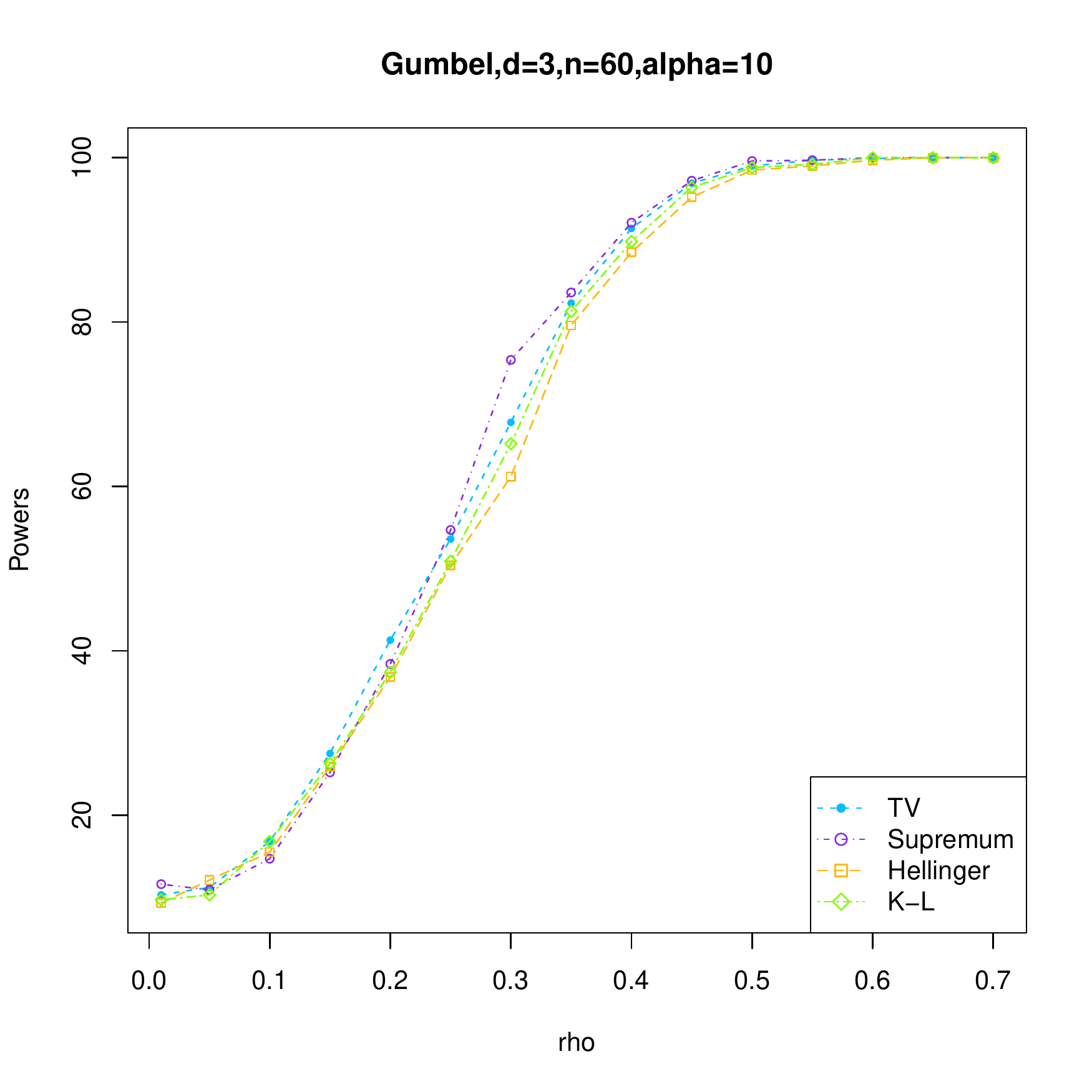}\hspace{.1cm}
\includegraphics[width=5.1cm,height=5.1cm]{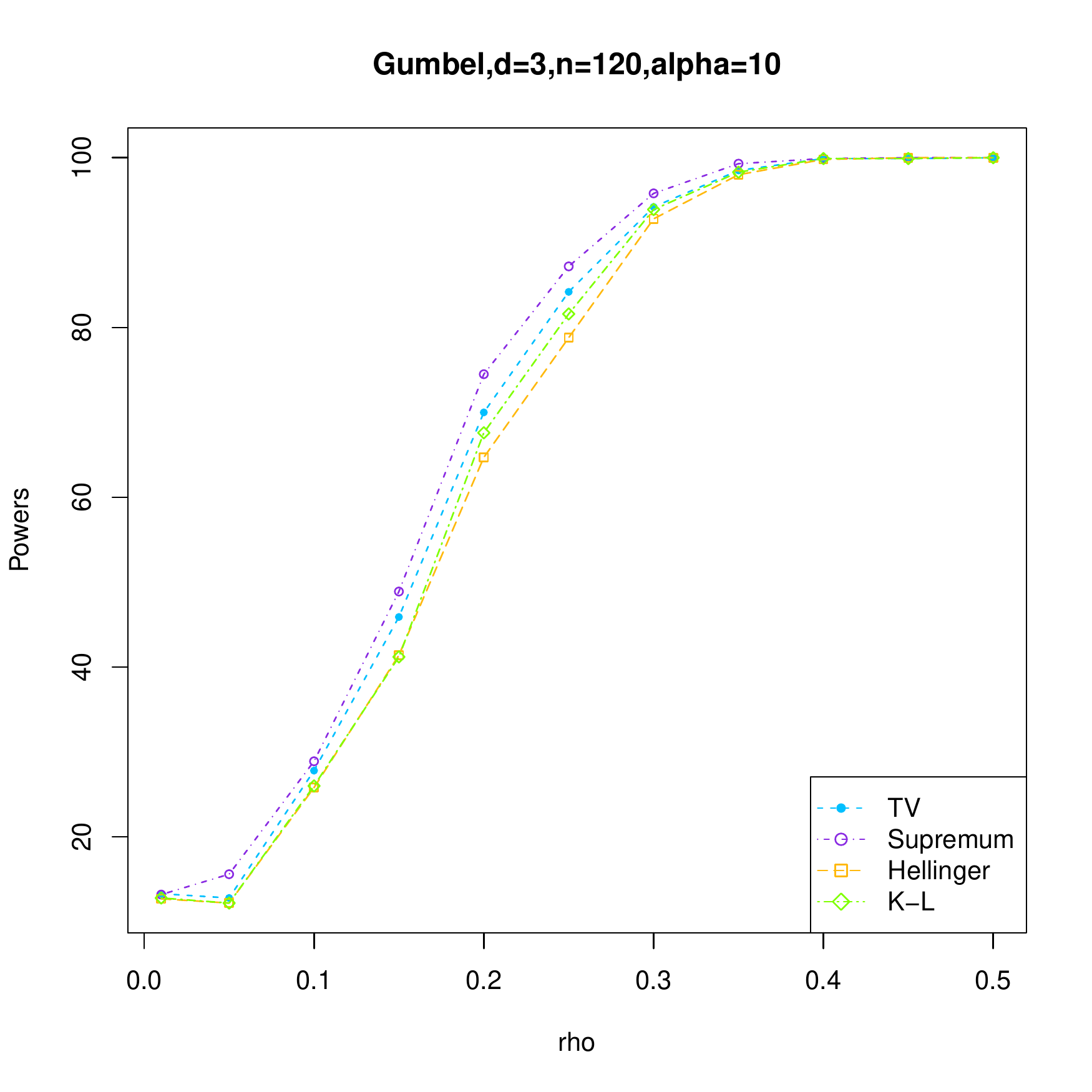}\hspace{.1cm}
\includegraphics[width=5.1cm,height=5.1cm]{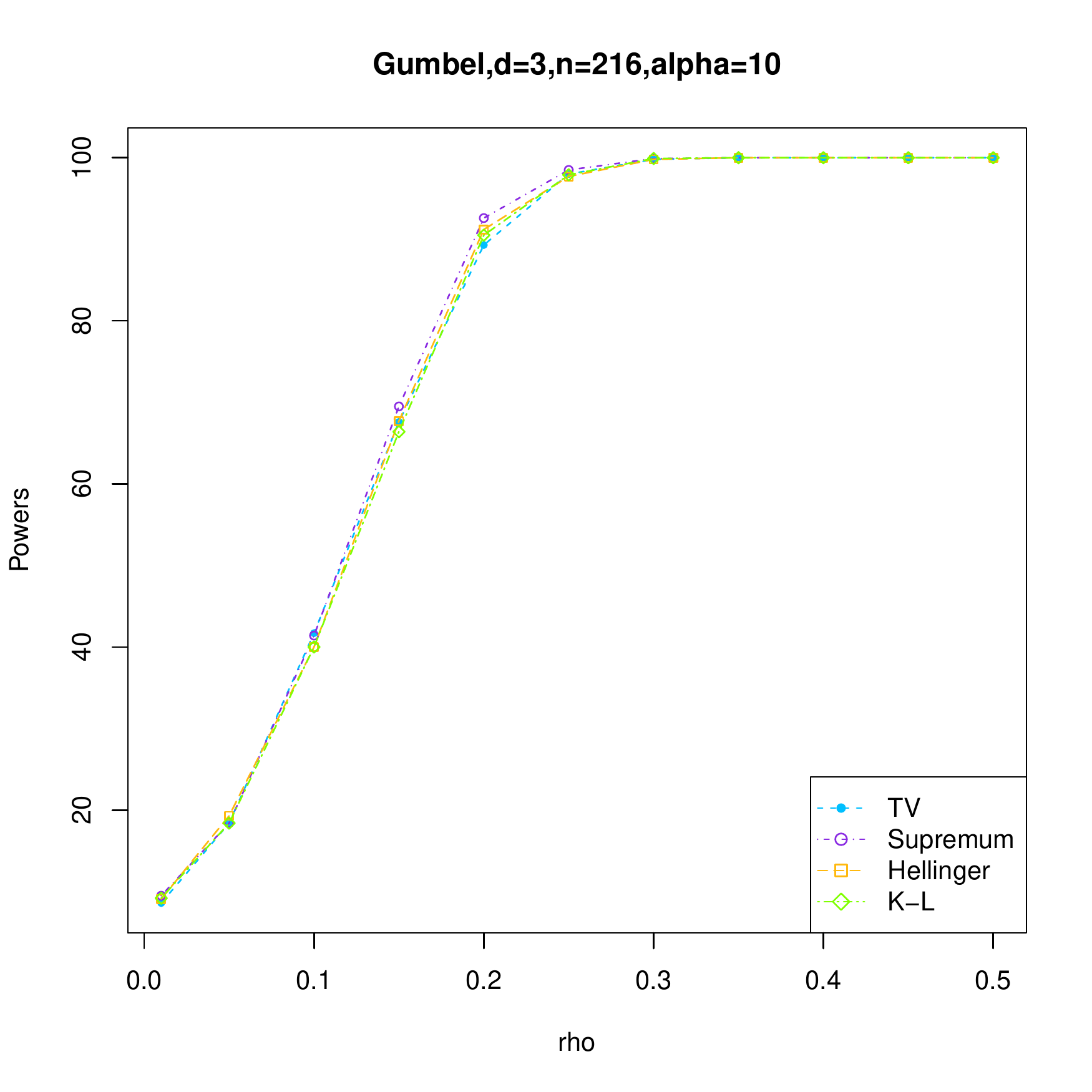}\hspace{.1cm}
\caption{Power comparison for the Gumbel family for various sample sizes; $d=3$.}
\label{Gumbel}
\end{center}
\end{figure}
\vspace*{-1ex}

\begin{figure}[h]
\begin{center}
\includegraphics[width=5.1cm,height=5.1cm]{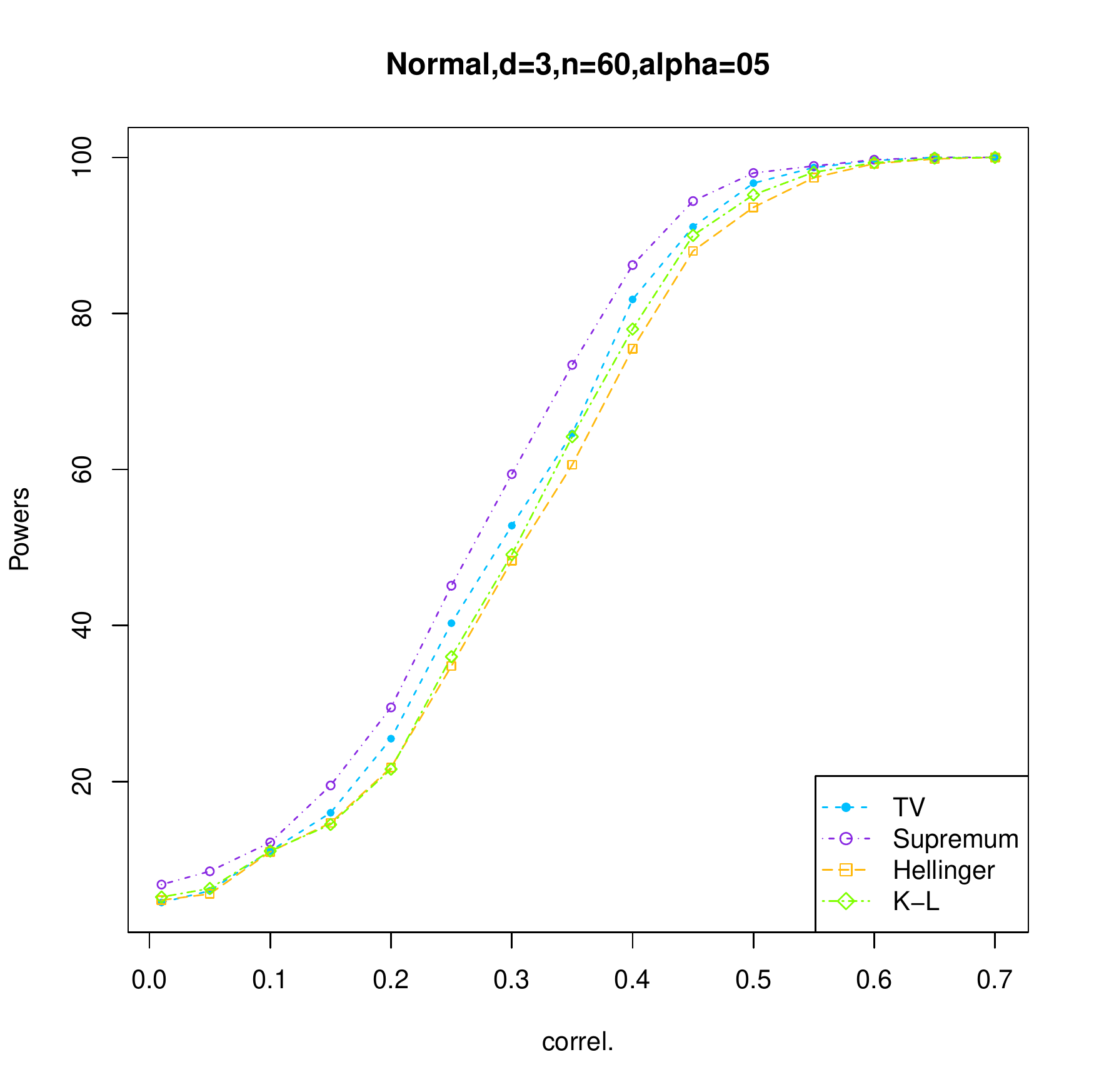}\hspace{.1cm}
\includegraphics[width=5.1cm,height=5.1cm]{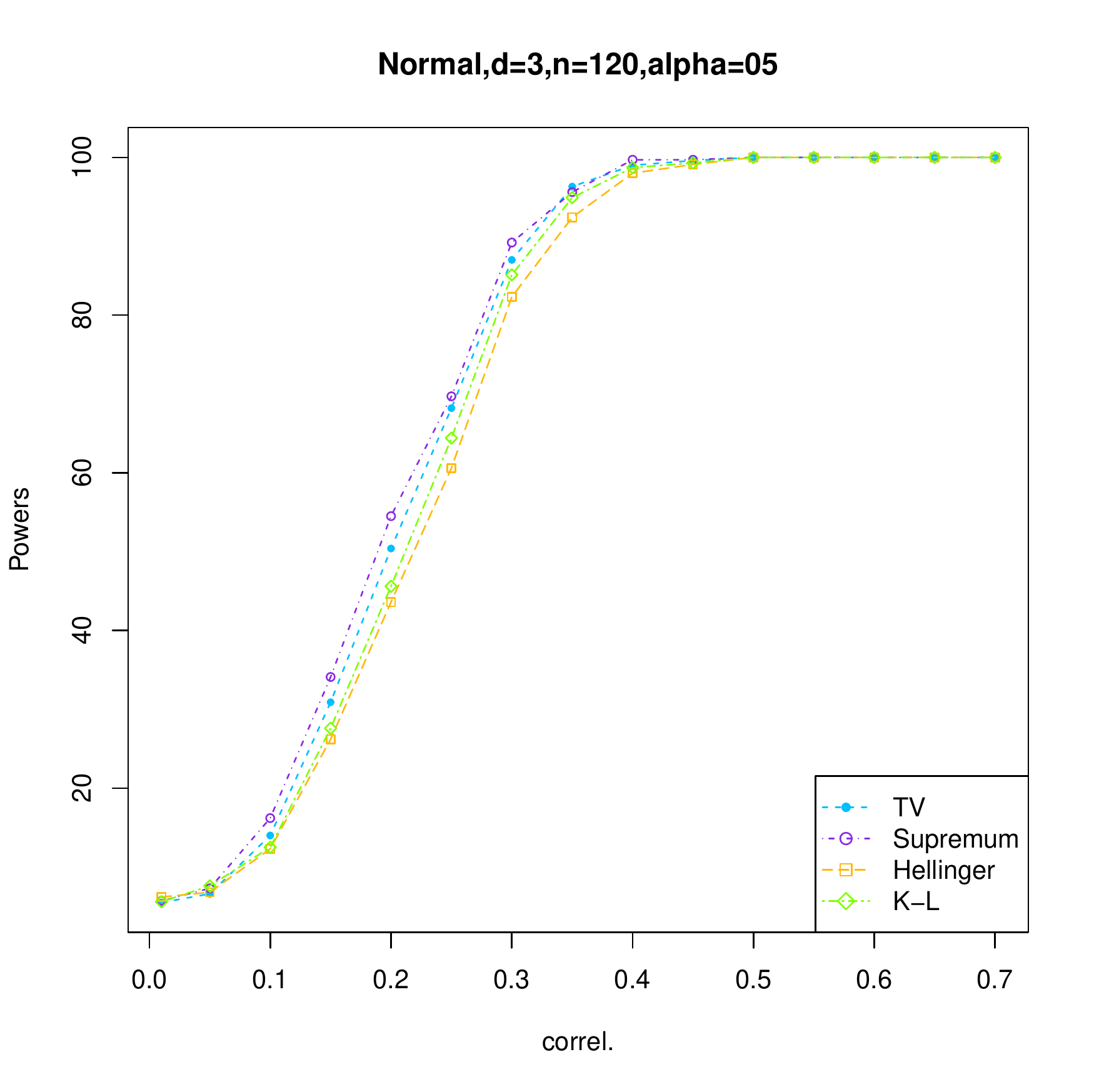}\hspace{.1cm}
\includegraphics[width=5.1cm,height=5.1cm]{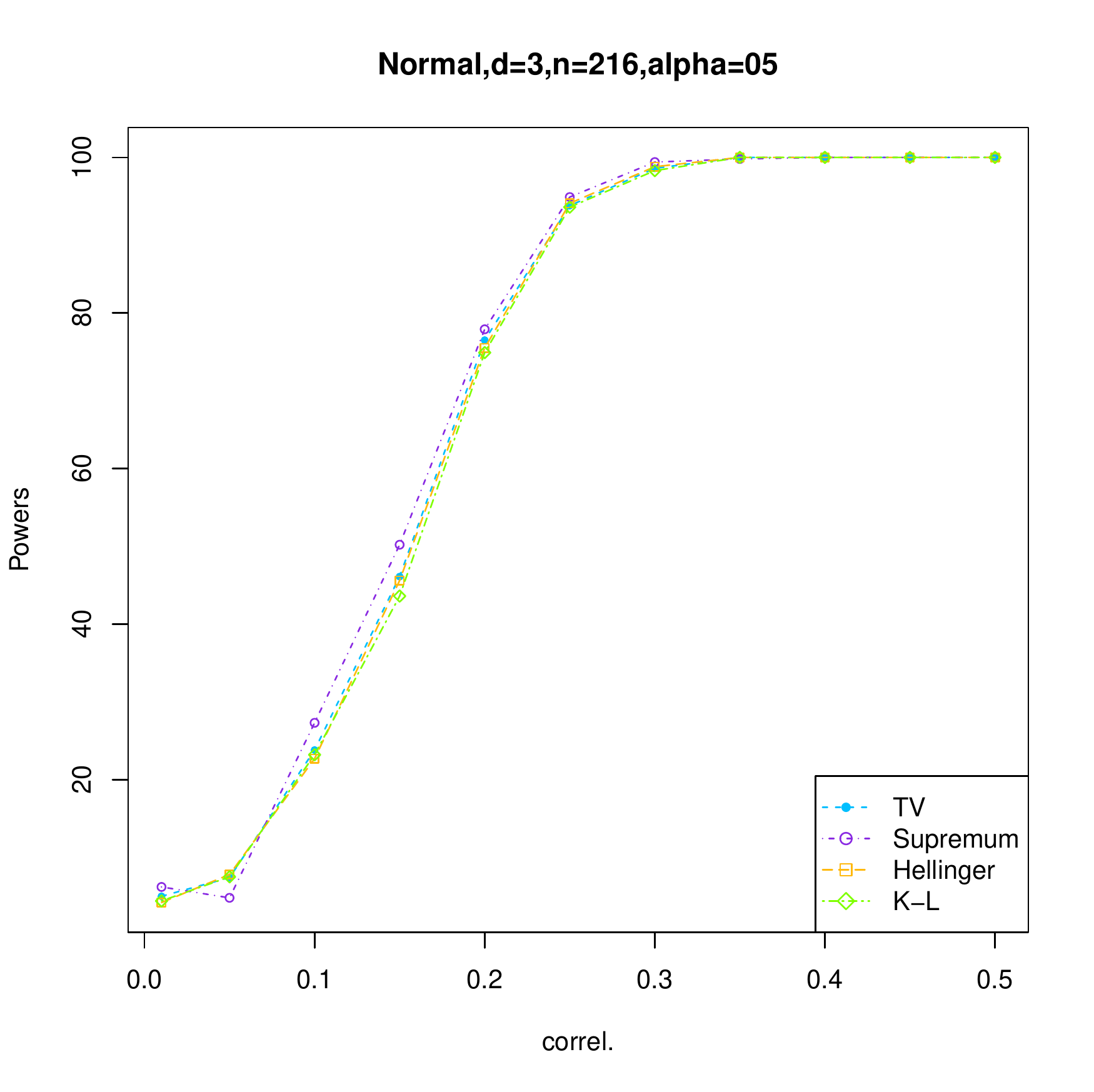}\hspace{.1cm}
\caption{Power comparisons for the Normal family for various sample sizes and pairwise correlations; $d=3$.}
\label{Normd3}
\end{center}
\end{figure}
\vspace*{-1ex}

\begin{figure}[h!]
\begin{center}
\includegraphics[width=5.1cm,height=5.1cm]{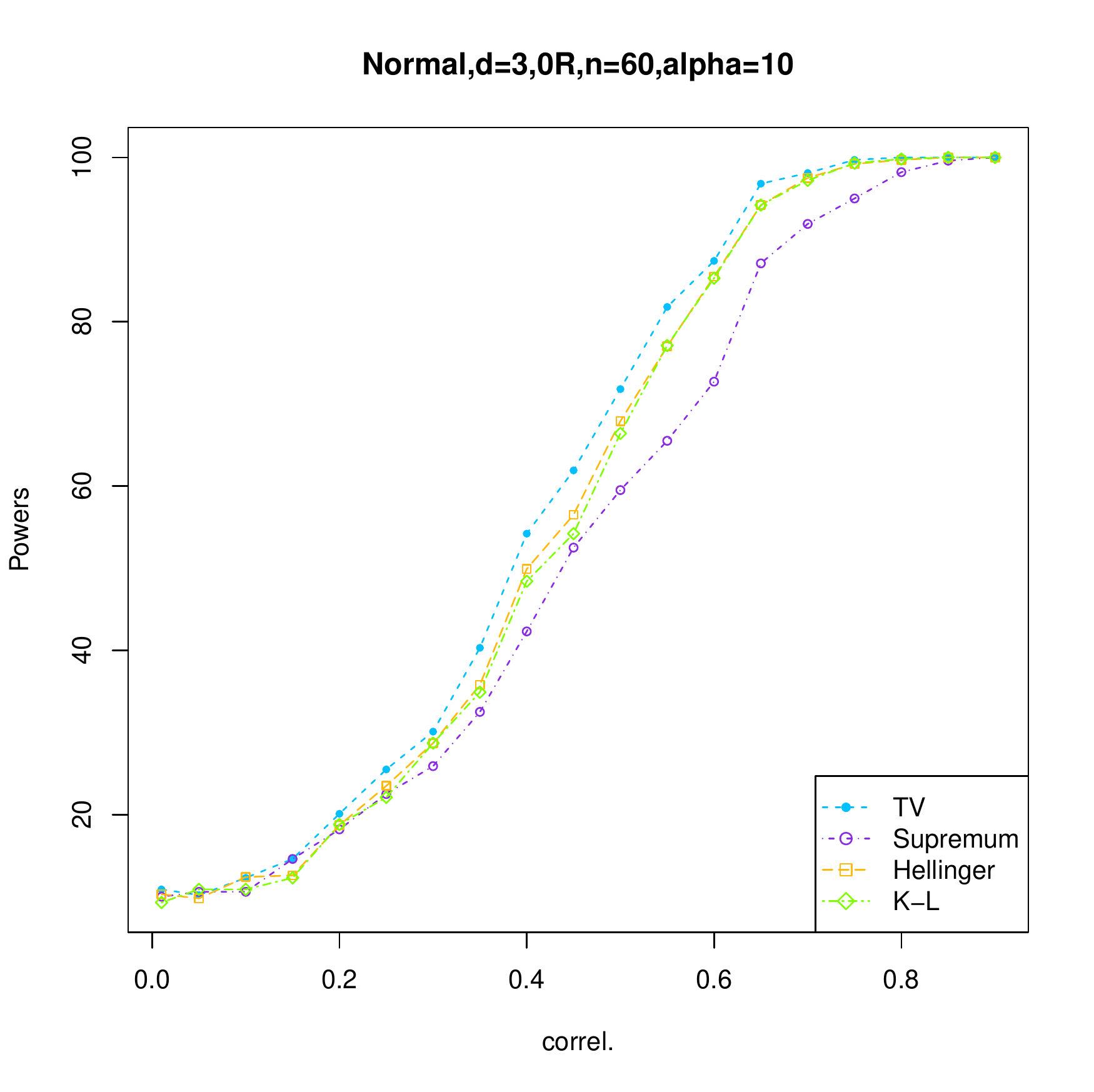}\hspace{.1cm}
\includegraphics[width=5.1cm,height=5.1cm]{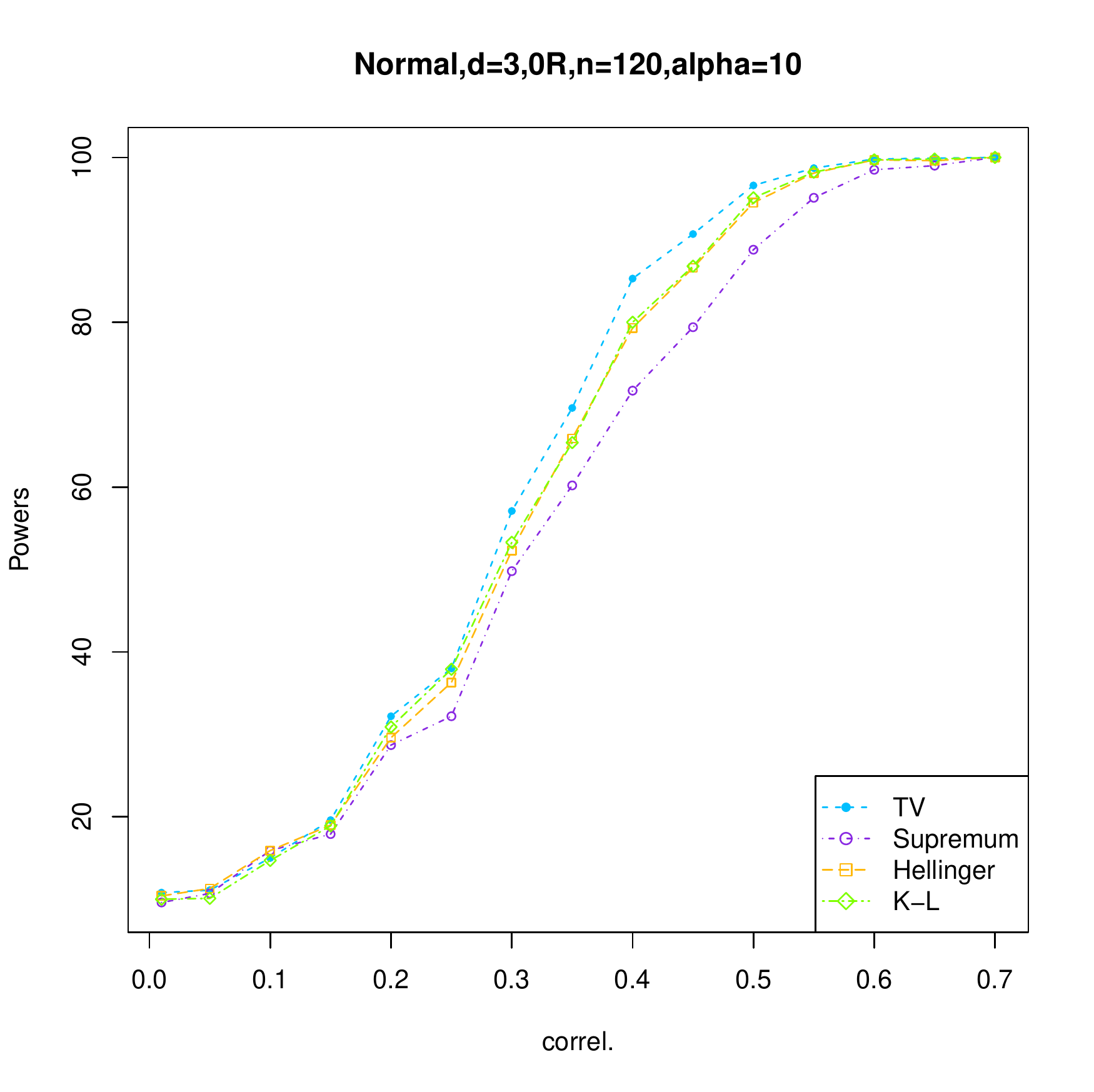}\hspace{.1cm}
\includegraphics[width=5.1cm,height=5.1cm]{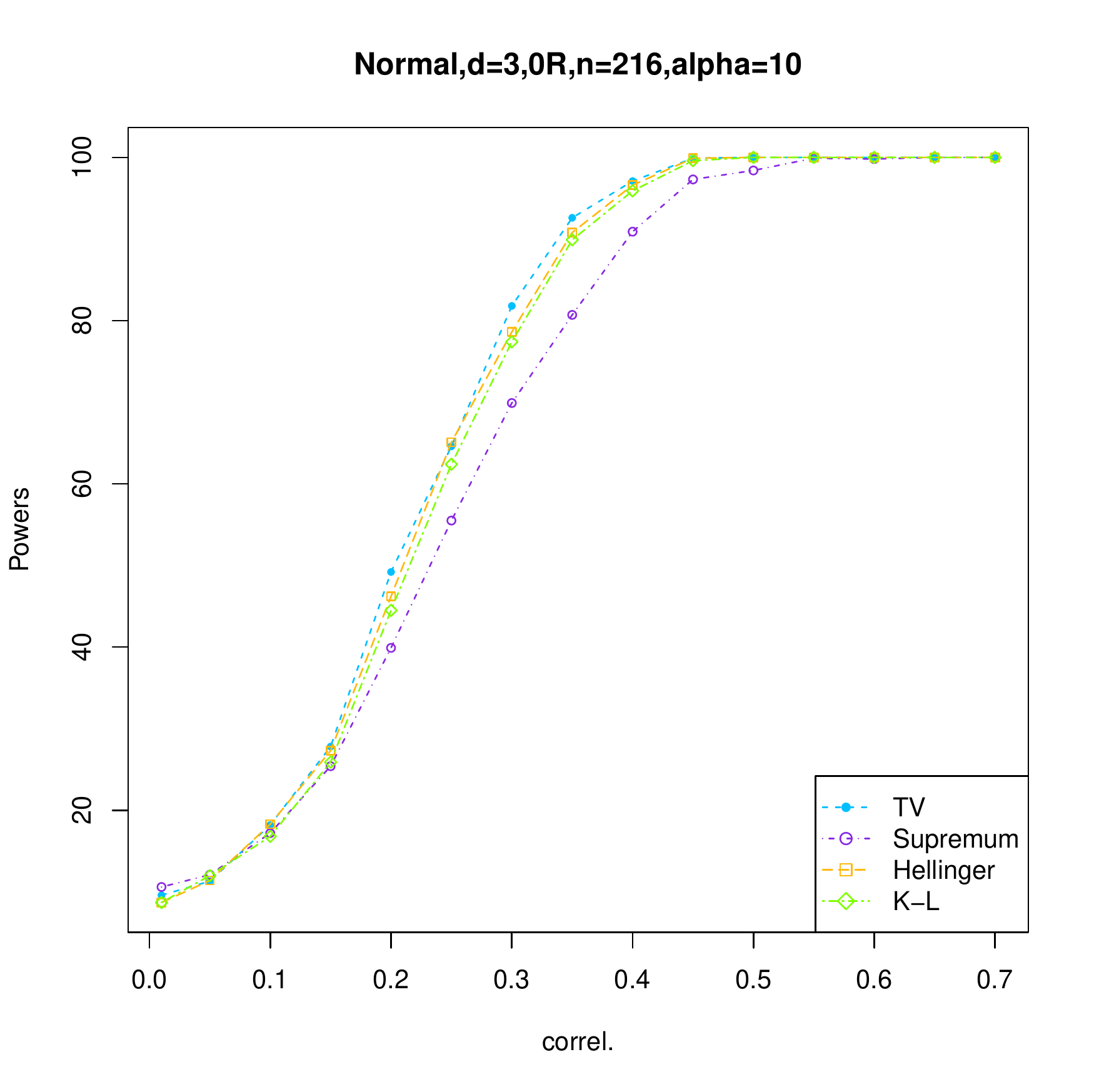}\hspace{.1cm}
\caption{Power comparisons for the Normal family for various sample sizes and pairwise correlations (weak dep.); $d=3$.}
\label{Normd30R}
\end{center}
\end{figure}

\begin{figure}[t!]
\begin{center}
\includegraphics[width=5.1cm,height=5.1cm]{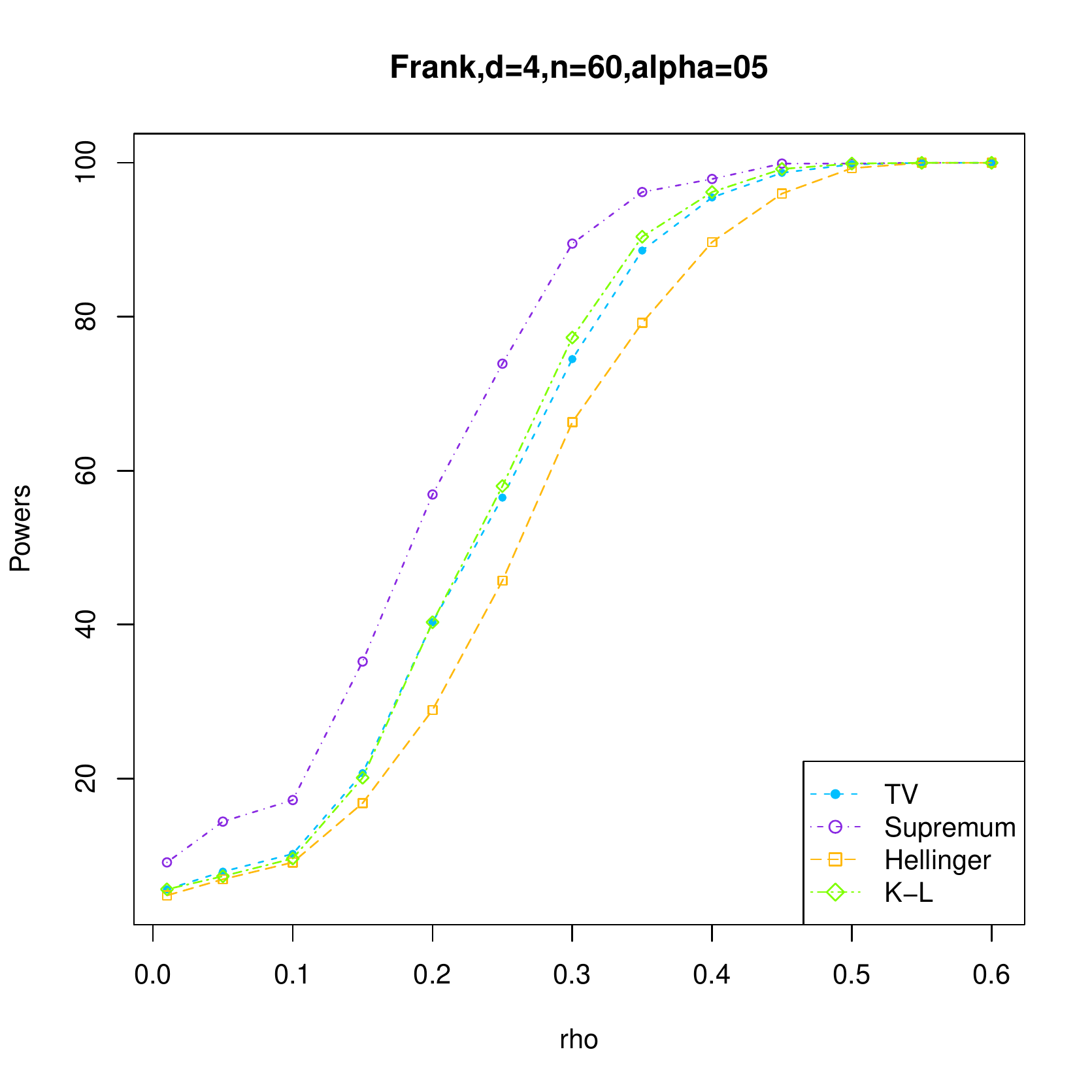}\hspace{.1cm}
\includegraphics[width=5.1cm,height=5.1cm]{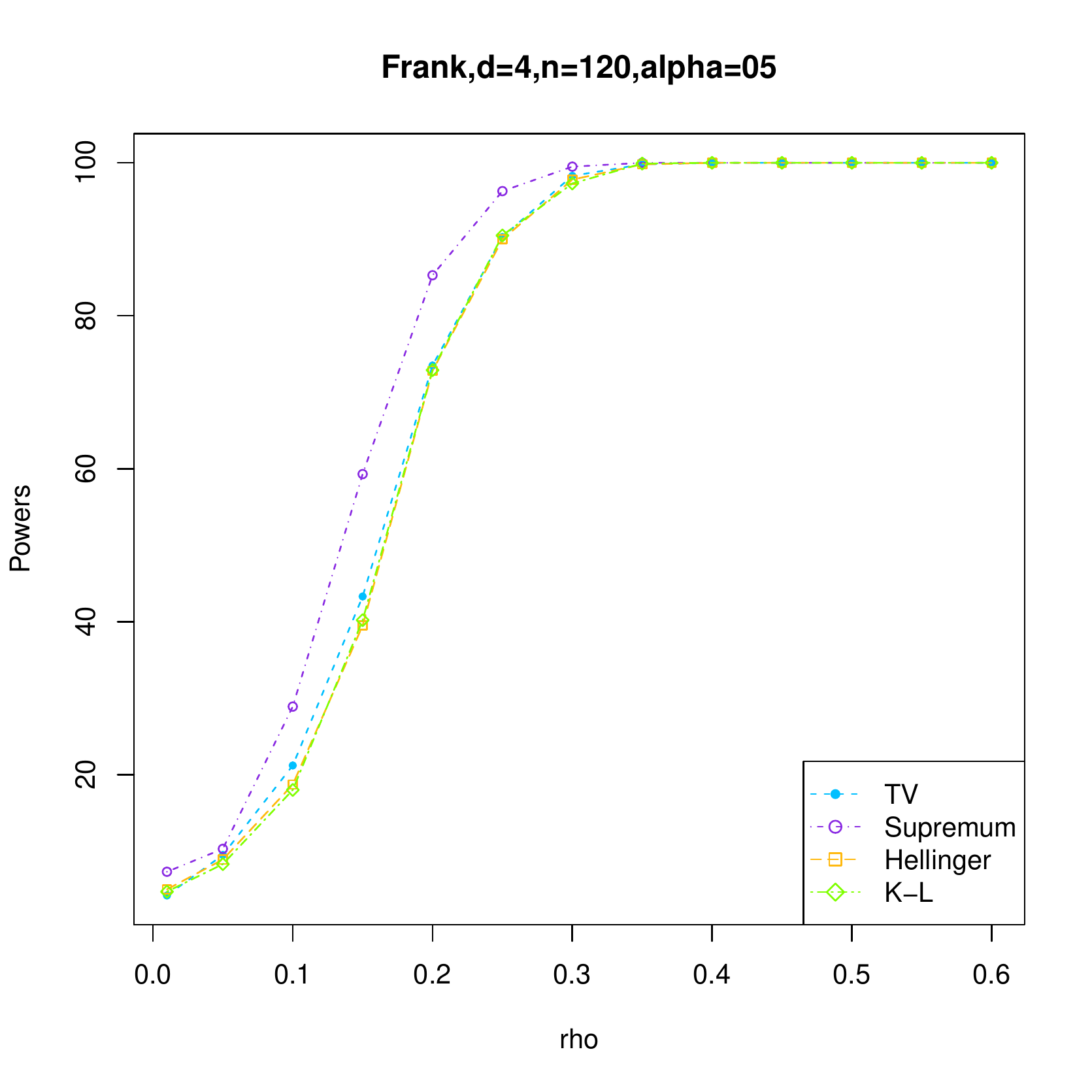}\hspace{.1cm}
\includegraphics[width=5.1cm,height=5.1cm]{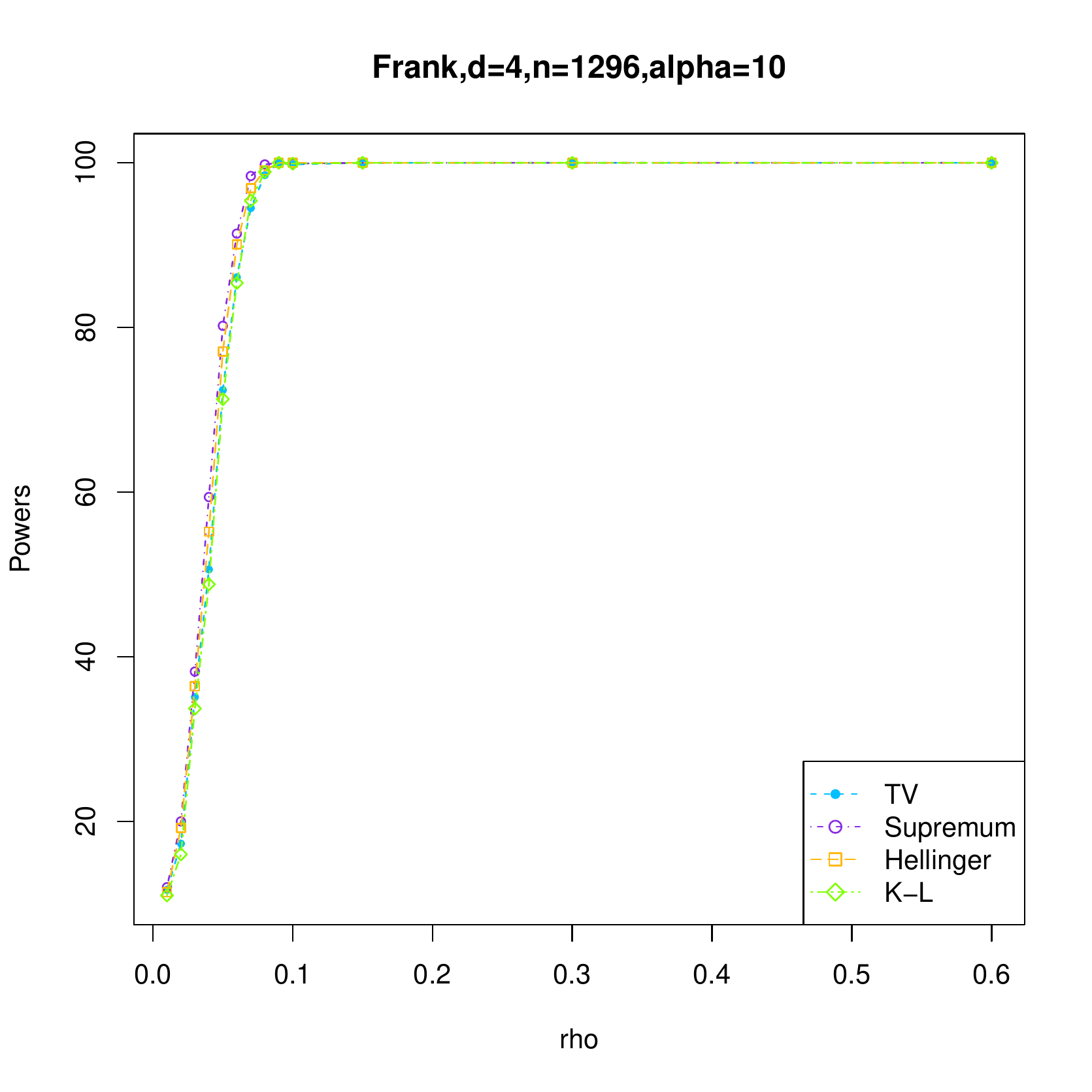}\hspace{.1cm}
\caption{Power comparisons for the Frank family for various sample sizes; $d=4$.}
\label{Frank}
\end{center}
\end{figure}
\vspace*{-1ex}

\begin{figure}[h]
\begin{center}
\includegraphics[width=5.1cm,height=5.1cm]{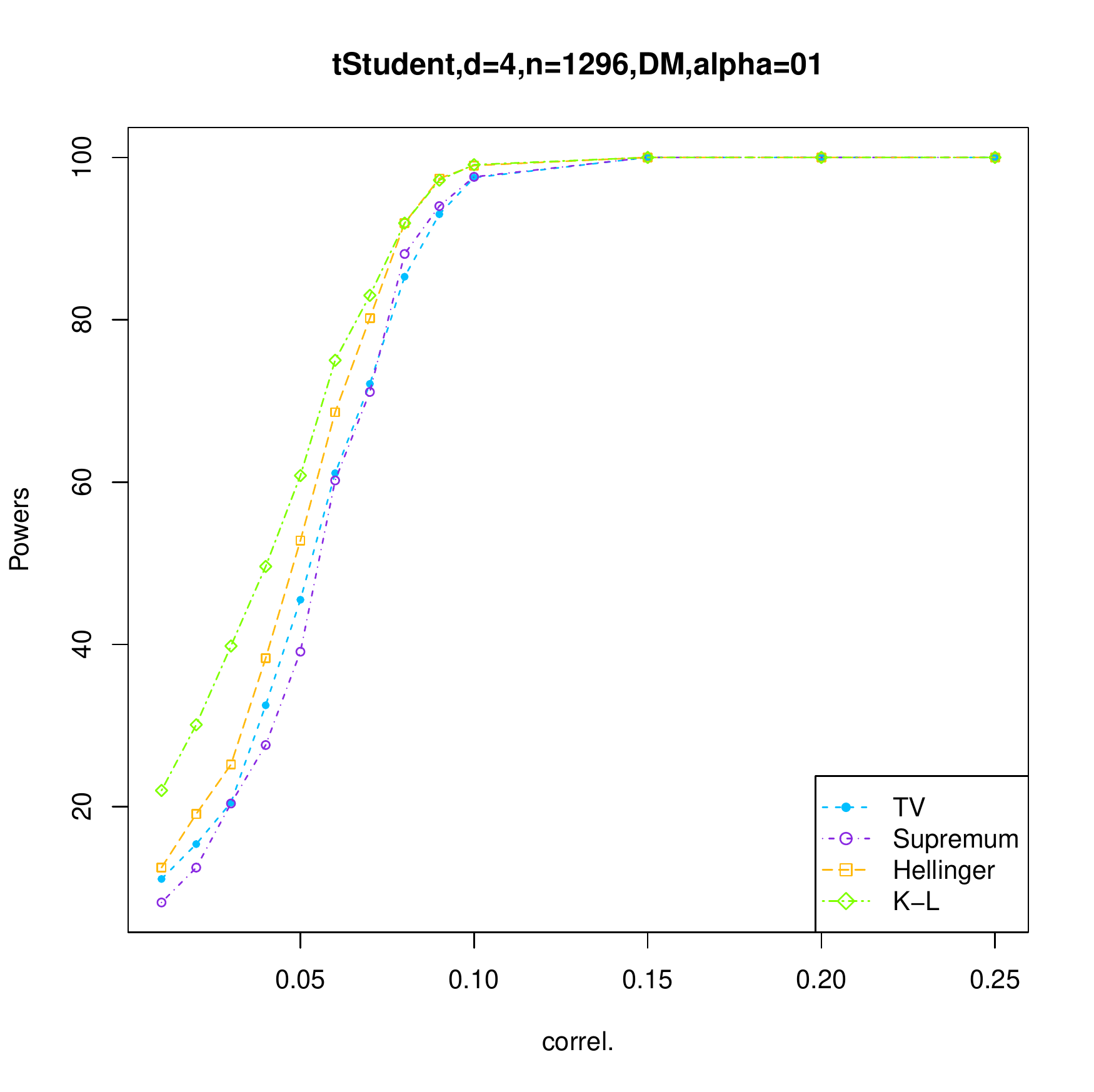}\hspace{.1cm}
\includegraphics[width=5.1cm,height=5.1cm]{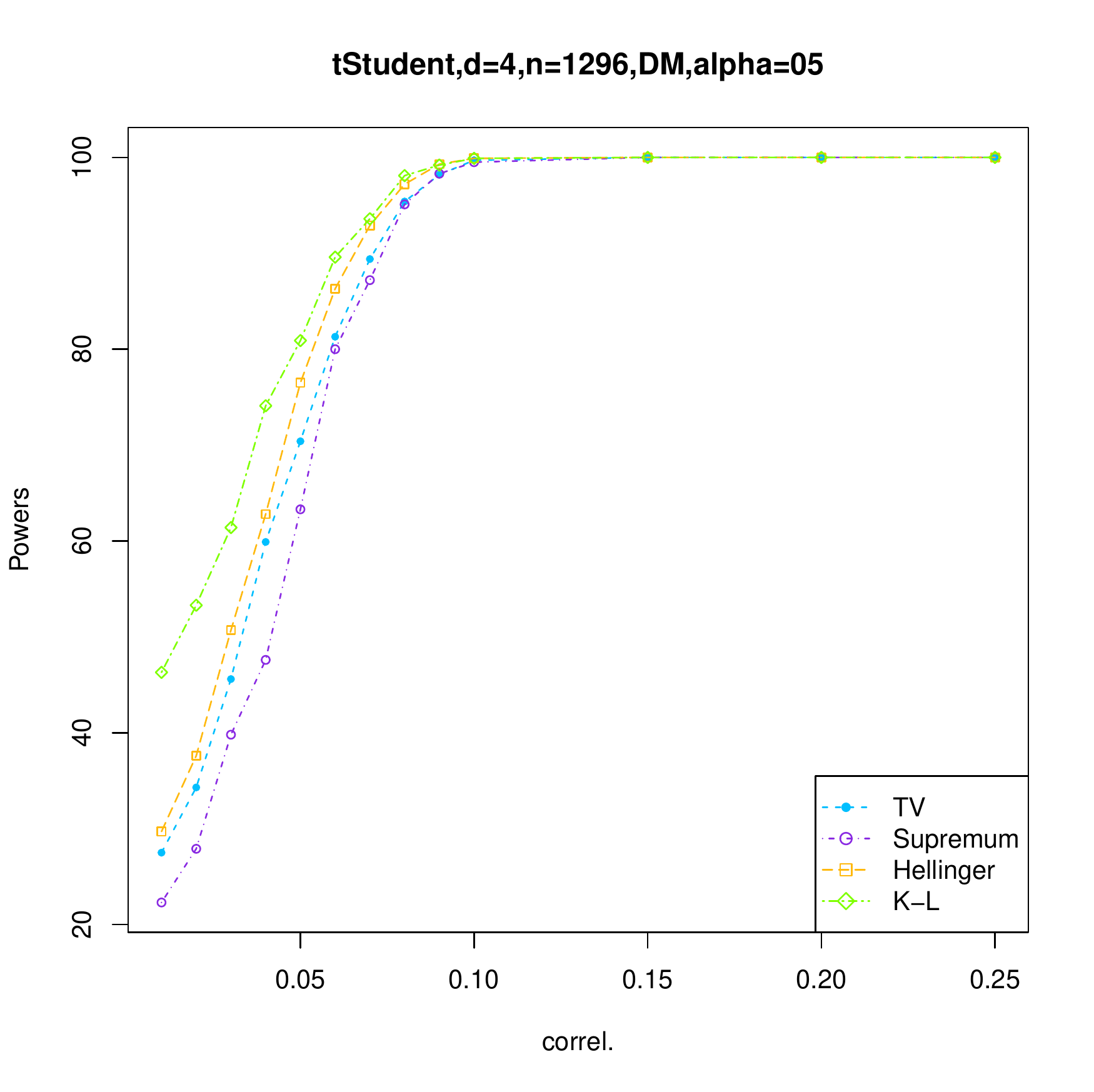}\hspace{.1cm}
\includegraphics[width=5.1cm,height=5.1cm]{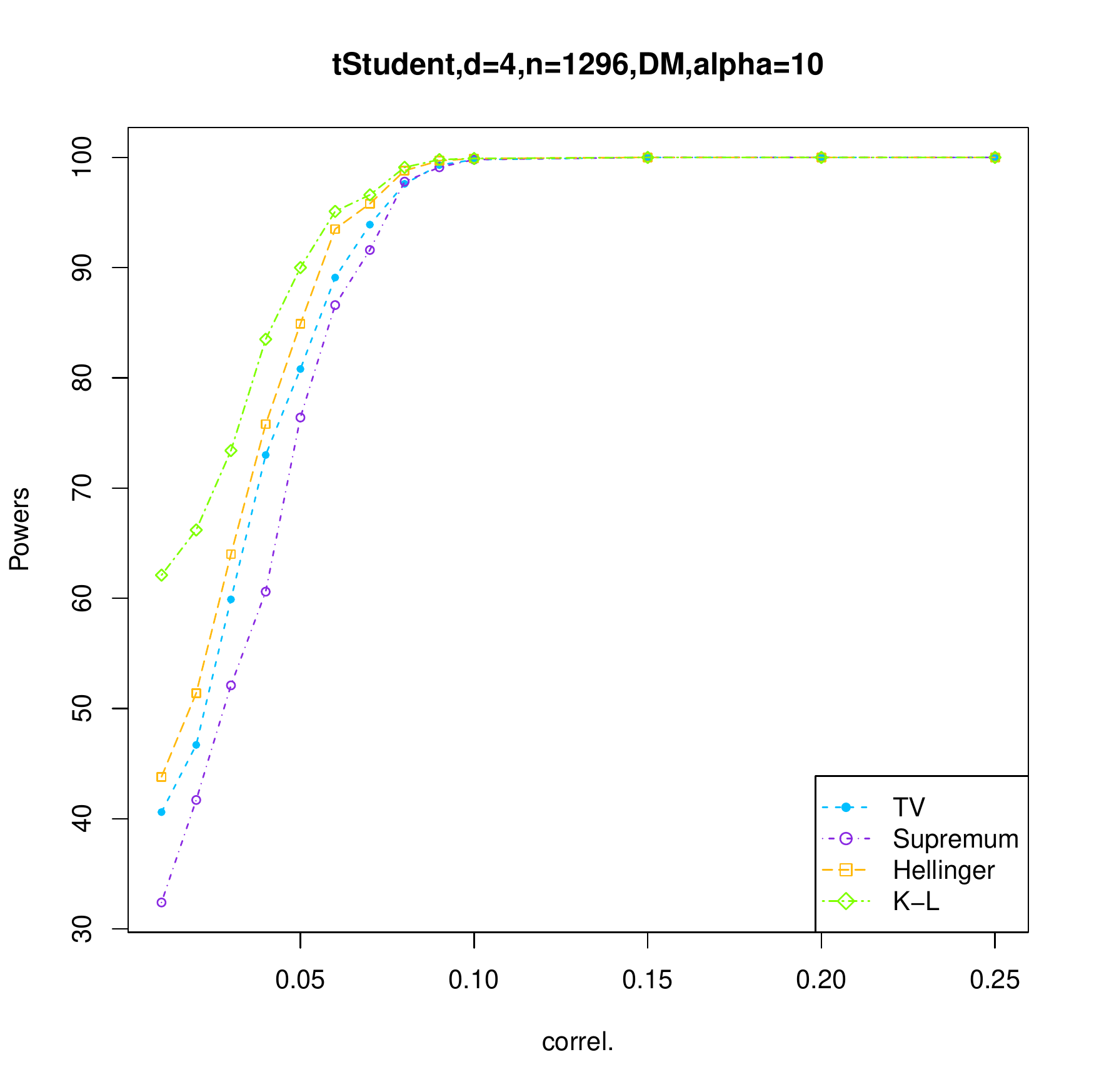}\hspace{.1cm}
\caption{Power comparisons for the Student $t$ family for various pairwise correlations; $n=1296$, $d=4$.}
\label{St4}
\end{center}
\end{figure}
\vspace*{-1ex}

\newpage

\end{document}